\documentclass{scrartcl}
\usepackage{amsmath}
\usepackage{amsfonts}
\usepackage{amssymb}
\usepackage{dsfont}
\usepackage[T1]{fontenc}
\usepackage[latin1]{inputenc}
\usepackage{epsfig}
\usepackage{graphicx}

\usepackage[boxed, lined]{algorithm2e}
\usepackage{float}
\usepackage{comment}

\newcommand{\ba}{\mathbf{a}}
\newcommand{\bb}{\mathbf{b}}
\newcommand{\bc}{\mathbf{c}}
\newcommand{\by}{\mathbf{y}}

\newcommand{\bB}{\mathbf{B}}
\newcommand{\bA}{\mathbf{A}}

\newcommand{\bID}{\mathbf{1}}
\newcommand{\D}{{\mathcal{D}}}

\newcommand{\Nu}{{\mathcal{N}}}

\newcommand{\N}{\mathbb{N}}

\newcommand{\R}{\mathbb{R}}

\newcommand{\Rd}{\mathbb{R}^d}

\newcommand{\beq}{\begin{eqnarray*}}

\newcommand{\eeq}{\end{eqnarray*}}

\newcommand{\beqm}{\begin{eqnarray}}

\newcommand{\eeqm}{\end{eqnarray}}

\newtheorem{theorem}{Theorem}

\newtheorem{lemma}{Lemma}
\newtheorem{definition}{Definition}

\newcommand{\EXP}{{\mathbf E}}
\newcommand{\PROB}{{\mathbf P}}

\renewcommand{\bf}{\normalfont \bfseries}
\renewcommand{\it}{\normalfont \itshape}

\allowdisplaybreaks
\begin{document}
\renewcommand{\thefootnote}{\fnsymbol{footnote}}
\newcommand{\F}{{\cal F}}
\newcommand{\Sp}{{\cal S}}
\newcommand{\G}{{\cal G}}

\begin{center}

  {\LARGE \bf
    On the rate of convergence of a neural network
    regression estimate learned by gradient descent
  }
\footnote{
Running title: {\it Neural network regression estimates}}
\vspace{0.5cm}
\footnote{This paper was presented at the Joint Statistical Meetings (JSM) 2019 in Denver, Colorado.}

Alina Braun$^{1}$, Michael Kohler$^{1,}$\footnote{Corresponding author. Tel: +49-6151-16-23382, Fax:+49-6151-16-23381}
and Harro Walk$^{2}$\\

{\it $^1$
Fachbereich Mathematik, Technische Universit\"at Darmstadt,
Schlossgartenstr. 7, 64289 Darmstadt, Germany,
email: braun@mathematik.tu-darmstadt.de, kohler@mathematik.tu-darmstadt.de}

{\it $^2$
Fachbereich Mathematik, Universit\"at Stuttgart,
Pfaffenwaldring 57, 70569 Stuttgart, Germany,
email: walk@mathematik.uni-stuttgart.de
}

\end{center}
\vspace{0.5cm}

\begin{center}
September 2, 2019
\end{center}
\vspace{0.5cm}

\noindent
    {\bf Abstract}\\
    Nonparametric regression with random design is considered.
    Estimates are defined by minimzing a penalized empirical $L_2$ risk
    over a suitably chosen class of neural networks
    with one hidden layer  via gradient
    descent. Here, the gradient descent procedure
    is repeated several times with
    randomly chosen starting values for the weights,
    and from the list of constructed
    estimates the one with the minimal empirical $L_2$ risk
    is chosen. Under the assumption that the
    number of randomly chosen starting values and the number of
    steps for gradient descent are sufficiently large it is shown
    that the resulting estimate achieves (up to a logarithmic factor)
    the optimal rate of convergence in a projection pursuit model.
    The final sample size performance of the estimates is illustrated
    by using simulated data. 

    \vspace*{0.2cm}


\vspace*{0.2cm}


\section{Introduction}
\label{se1}
\subsection{Scope of this article}
Motivated by the huge success of multilayer neural networks
in applications 
(see, e.g.,
Schmidhuber (2015) and the literature cited therein)
there has been an increasing interest in the
theoretical analysis of such estimates. Often this is done
in the area of nonparametric regression, and recently there
has been a tremendous progress in the theoretical understanding
of least squares regression estimates based on deep neural networks, i.e.,
neural networks with many hidden layers. The corresponding theoretical
results are based on the derivation of new approximation results for piecewise
polynomials by neural networks, and they make extensive use of the network
structure, which allows to exploit compository assumptions on the
structure of the regression function in order to circumvent
the curse of dimensionality (cf.,
Kohler and Krzy\.zak (2017),
Bauer and Kohler (2017),
Schmidt-Hieber (2017),
Imaizumi and  Fukumizu (2018),
Eckle and  Schmidt-Hieber (2018) and
Kohler, Krzy\.zak and Langer (2019)).

In all the articles above the neural network regression estimate is defined
as a nonlinear least squares estimate, i.e., as a function which minimizes
the empirical $L_2$ risk
over a nonlinear class of neural
networks. In practice, it is usually not possible to find the
global minimum of the empirical $L_2$ risk over a
nonlinear class of neural networks and 
one usually tries to find a local minimum using, for instance, the
gradient descent algorithm. So although the above theoretical
results are quite impressive, there is a big gap between the
estimates studied theoretically and the estimates used in
practice. 

The purpose of this paper is to narrow this gap.
To do this, we consider the following question:
If we define a neural network regression estimate theoretically
exactly as it is implemented in practice, can we show a
rate of convergence result for this estimate?
The ultimative goal is to analyze theoretically neural network regression
estimates which are actually used in practice. As a first step
in this direction we define a simple neural network regression
estimate where we use gradient descent in order to learn the
weights of a neural network with one hidden layer in
a projection pursuit model. We show that if we repeatedly
apply this procedure to starting values, which are chosen randomly from
a special structure, then, for sufficiently many starting values
and steps in each procedure, we will find an estimate which
achieves the optimal rate
of convergence up to a logarithmic factor in  this projection pursuit model.

\subsection{Nonparametric regression}

We study neural network estimates
in the context of nonparametric regression with random
design. Here,
$(X,Y)$ is an $\Rd \times \R$--valued random vector
satisfying $\EXP \{Y^2\}<\infty$, and given a sample
of $(X,Y)$ of size $n$, i.e., given a data set
\begin{equation}
  \label{se1eq1}
\D_n = \left\{
(X_1,Y_1), \dots, (X_n,Y_n)
\right\},
\end{equation}
where
$(X,Y)$, $(X_1,Y_1)$, \dots, $(X_n,Y_n)$ are i.i.d. random variables,
the aim is to construct an estimate
\[
m_n(\cdot)=m_n(\cdot, \D_n):\Rd \rightarrow \R
\]
of the regression function $m:\Rd \rightarrow \R$,
$m(x)=\EXP\{Y|X=x\}$ such that the $L_2$ error
\[
\int |m_n(x)-m(x)|^2 \PROB_X (dx)
\]
is ``small'' (see, e.g., Gy\"orfi et al. (2002)
for a systematic introduction to nonparametric regression and
a motivation for the $L_2$ error).

It is well--known that
one
needs smoothness assumptions on the regression function in
order to derive non--trivial results
on
the rate of convergence
of nonparametric regression estimates
(cf., e.g., Theorem 7.2 and Problem 7.2 in
Devroye, Gy\"orfi and Lugosi (1996) and
Section 3 in Devroye and Wagner (1980)).
To do this we will use the following definition.
\begin{definition}
\label{intde2}
  Let $p=q+s$ for some $q \in \N_0$ and $0< s \leq 1$,
where $\N_0$ is the set of nonnegative integers.
A  function $f:\R^d \rightarrow \R$ is called
\textbf{$(p,C)$-smooth}, if for every $\alpha=(\alpha_1, \dots, \alpha_d) \in
\N_0^d$
with $\sum_{j=1}^d \alpha_j = q$ the partial derivative
$\frac{
\partial^q f
}{
\partial x_1^{\alpha_1}
\dots
\partial x_d^{\alpha_d}
}$
exists and satisfies
\[
\left|
\frac{
\partial^q f
}{
\partial x_1^{\alpha_1}
\dots
\partial x_d^{\alpha_d}
}
(x)
-
\frac{
\partial^q f
}{
\partial x_1^{\alpha_1}
\dots
\partial x_d^{\alpha_d}
}
(z)
\right|
\leq
C
\cdot
\| x-z \|^s
\]
for all $x,z \in \R^d$, where $\Vert\cdot\Vert$ denotes the Euclidean norm.
\end{definition}
Stone (1982) showed that the optimal minimax rate of convergence in nonparametric
regression for $(p,C)$-smooth functions is $n^{-2p/(2p+d)}$. In case
that $d$ is large compared to $p$ this rate of convergence is
rather slow (so called curse of dimensionality). In the sequel we want
to
circumvent this curse of dimensionality by imposing the additional
constraint on the regression function that it satisfies a projection
pursuit model, i.e., by assuming that it satisfies
\begin{equation}
\label{se1eq2}
m(x)=\sum_{s=1}^r g_s(\bc_s^T x) \quad (x \in \Rd)
\end{equation}
for some $r \in \N$, $\bc_s \in \Rd$, where $\|\bc_s\| = 1$, and
$(p,C)$-smooth functions
$g_s:\R \rightarrow \R$ $(s=1,
\dots,r)$. Under this assumption our aim is to show that suitably
defined neural network estimates, which can be actually implemented
in an application, can achieve the one-dimensional rate of convergence.

\subsection{Main result of this article}
In this paper we study neural network regression estimates
using neural networks with one hidden layer in the above projection
pursuit model, i.e., we assume that the regression function
satisfies (\ref{se1eq2}). We learn the weights of our neural
network regression estimate by choosing in a first step
randomly vectors
for the directions $\bc_s$ of our projection pursuit model,
by defining in a second step an appropriate starting value for
the weights of our neural network regression estimate based
on the randomly chosen directions, and by applying in a third
step successively many gradient descent steps in order to optimze the
weights of our neural network. Then we repeat this whole
procedure several times and choose from the list of 
estimates which we get the one with the minimal error
on our training data. 

Our main result is that for a sufficiently
large number of repititions of this procedure and a sufficiently
large number of gradient descent steps
the expected $L_2$ error of a  truncated version
of our  estimate converges towards zero
 in the projection pursuit
model (\ref{se1eq2}) in case of $(p,C)$--smoth functions
$g_s$ (where $p \leq 1$) with the rate of convergence
\[
\left(
\frac{(\log n)^3}{n}
\right)^{\frac{2p}{2p+1}},
\]
i.e., with the optimal rate
of convergence up to a logarithmic factor. 
Here, the rate of convergence is independent
of the dimension $d$ of $X$. Hence, our neural network regression
estimate is able to circumvent the curse of dimensionality in
the projection pursuit model (\ref{se1eq2}). 

We achieve this result by choosing our initial weights such that
the initial network basically  computes a piecewise constant
function and by showing that in this case the gradient descent
is able to choose the outer weights in the neural network
in an optimal way (provided the number of gradient descent
steps is sufficiently large).

\subsection{Discussion of related results}
It is well-known that it is possible
to circumvent the curse of dimensionality by
imposing additional constraints on the structure of the regression
function. Stone (1985) assumed that the regression function is
additive, i.e., that $m:\Rd \rightarrow \R$ satisfies
\[
m(x^{(1)},\dots,x^{(d)})
=
m_1(x^{(1)})+ \dots + m_d(x^{(d)})
\quad (x^{(1)},\dots,x^{(d)} \in \R)
\]
for some $(p,C)$--smooth univariate functions $m_1, \dots, m_d:\R
\rightarrow \R$, and showed that in this case suitably defined
spline estimates achieve the corresponding univariate rate of
convergence.
Stone (1994) extended this results to interaction models, where
the regression function is assumed to be a sum of functions
applied to at most $d^* < d$ components of $x$ and showed that
in this case suitably defined spline estimates
achieve
the $d^*$--dimensional rate of convergence.
Other classes of functions which enable us to achieve a better rate of convergence
results include single index models, where
\[
m(x)=g(\bc^T x) \quad (x \in \Rd)
\]
for some $\bc \in \Rd$ and $g:\R \rightarrow \R$
 (cf., e.g.,  H\"ardle and Stoker (1989), H\"ardle, Hall and Ichimura (1993), 
 Yu and Ruppert (2002), Kong and Xia (2007)  and  Lepski and Serdyukova (2014)),
and projection pursuit, where it is assumed that (\ref{se1eq2}) holds
for some $r \in \N$, $\bc_s \in \Rd$ and $g_s:\R \rightarrow \R$ $(s=1,
\dots,r)$ (cf., e.g., Friedman and Stuetzle (1981) and Huber (1985)).
Horowitz and Mammen (2007) studied the case of a regression function, which satisfies
\[
m(x)=g\left(\sum_{l_1=1}^{L_1}g_{l_1}  \left(\sum_{l_2=1}^{L_2}g_{l_1,l_2}\left( \ldots \sum_{l_r=1}^{L_r}g_{l_1,\ldots, l_r}(x^{l_1,...,l_r}) \right)\right)\right),
\]
where $g, g_{l_1}, \ldots, g_{l_1,\ldots, l_r}$ are $(p,C)$-smooth
univariate functions and $x^{l_1,...,l_r}$ are single components of
$x\in\Rd$ (not necessarily different for two different indices
$(l_1,\ldots,l_r)$). With the use of a penalized least squares
estimate for smoothing splines, they proved the rate $n^{-2p/(2p+1)}$.

For the $L_2$ error of a
single hidden layer neural network, Barron (1993, 1994) proved
 the dimensionless rate of convergence $n^{-1/2}$
(up to some logarithmic factor), provided the Fourier transform has a finite first
moment (which basically
requires that the function becomes smoother with increasing
dimension $d$ of $X$).
Restricting their study to the use of a certain cosine squasher as the activation function,
 McCaffrey and Gallant (1994) showed a rate of $n^{-\frac{2p}{2p+d+5}+\varepsilon}$ for the $L_2$  error of suitably defined single hidden layer neural network estimate for $(p,C)$-smooth functions.

Recently it was shown in several papers that neural networks
can achieve a dimensionality reduction in case that the regression
function is a composition of (sums of)
functions, where each of the function
is a function of at most $d^*<d$ variables.
The first paper in this respect was
Kohler and Krzy\.zak (2017), where it was shown that
under a corresponding assumption
 suitably defined multilayer
neural networks achieve the rate of convergence  $n^{-2p/(2p+d^*)}$
(up to some logarithmic factor) in case $p \leq 1$.
Bauer and Kohler (2017) showed that this result even holds
for $p>1$ provided the squashing function is suitably
chosen. Schmidt-Hieber (2017) showed similar results for neural
networks with ReLU activation function.
Eckle and Schmidt-Hieber (2018) showed that neural networks
with ReLU activation function can approximate well piecewise
polynomials with rather general partitions based on the intersection
of hyperplanes and used this result to relate the error of neural
network
estimates to the error of piecewise polynomial partitioning estimates.
Kohler, Krzy\.zak and Langer (2019) derived a similar result for
neural networks with squashing functions as activation function
and used this result to prove that neural networks
are able to circumvent the curse of dimensionality in case that
the regression function has a low local dimensionality.
Results concerning the approximation of piecewise polynomials
with partitions with rather general smooth boundaries
by neural networks
 have been
derived in Imaizumi and Fukamizu (2018).

The above mentioned results show that least squares neural network
regression estimates are able to circumvent the curse of
dimensionality
under much more general assumptions than the projection
pursuit model assumed in this paper. However, these
estimates cannot be computed in practice, whereas our result
shows that in the projection pursuit model we can
achieve this with neural networks even in the case where
we restrict ourselves to estimates which can be 
computed much easier.

Gradient descent has been studied in many different papers,
see, e.g., Karimi, Nutini and Schmidt (2018) and the literature
cited therein.
A standard reference is the monograph Luenberger and Ye (2016). We also
mention Poljak (1981) as an early paper, where the case of noise corrupted
function values is considered, too. Stochastic approximation deals with the
latter field, see, e.g., the monograph Kushner and Yin (2003), and here in a
classic situation the constant factor at the gradient is replaced by a
decreasing factor at a vector of divided differences (multidimensional
Kiefer-Wolfowitz method). The paper of White (1989, 1992) brings together
the two fields of stochastic approximation and neural network models (see
also Fabian (1994)).
In Dippon and Fabian (1994) and Dippon (1998)
it is explained how gradient descent
in stochastic approximation
can be combined with a slowly convergent
global optimizer in order to find not only a local but even
a global minimum of a general function. The main difficulty
of using such results to derive rate of convergence results
for neural network regression estimates lies in the fact that for 
neural network regression estimates
the neural network is using more and more neurons with increasing
sample size. This means that it is not sufficient to analyze gradient
descent applied to a fixed function where the number of steps 
is tending to
infinity. Instead the function is changing for increasing number of
steps.
Basically, this requires the ability to analyze the behaviour of
gradient descent for a finite number of steps. As far as we know
such results do not exist in the literature in case of a general
function like the empirical $L_2$ risk of a neural network
(which is neither convex nor has a global minimum or an easily
analysable Hessian matrix considered as a real-valued function
of the weight vector).

There are quite a few articles in computer science
where people try to prove that
backpropagation leads to good neural network estimates.
Unfortunately, the approaches used there do not lead to similarly powerful
rate of convergence results for neural networks
as in our article here.
For instance, 
Arora et al. (2018), 
Kawaguchi (2016),
and Du and Lee (2018)
analyzed gradient descent for neural networks with linear or quadratic activation function.
For such neural networks there do not exist good approximation results, 
consequently, one cannot derive from these results 
a rate of convergence result 
comparable 
to that in our article.
Du et al. (2018) 
analyzed gradient descent applied to neural networks
with one hidden layer in case of a Gaussian input distribution. 
They used the expected gradient instead of the gradient 
in their gradient descent routine, 
and therefore, 
their result cannot be used to derive a rate of convergence result
for an estimate learned by gradient descent as the one in our paper.
Liang et al. (2018)
applied gradient descent to a modified loss function in classification, 
where it is assumed that the data can be interpolated by a neural network.
Here, the second assumption is not satisfied in nonparametric regression
and it is unclear whether the main idea (of simplifying the estimation by a modification of the loss function) 
can also be used in a regression setting.
In Allen-Zhu, Li and Song (2019),
also Kawaguchi and Huang (2019),
it is shown that gradient descent leads to 
a small empirical $L_2$ risk in overparametrized neural networks. 
Here, it is unclear what the $L_2$ risk of the estimate is (and a bound on this term is necessary in order to derive results like in our paper).
In particular, due to the fact that the networks are overparametrized,
a bound on the empirical $L_2$ risk might be not useful for bounding the $L_2$ risk.

\subsection{Notation}
Throughout the paper, the following notation is used:
The sets of natural numbers, natural numbers including $0$ and real numbers
are denoted by $\N$, $\N_0$ and $\R$, respectively. For $z \in \R$, we denote
the smallest integer greater than or equal to $z$ by
$\lceil z \rceil$.
The Euclidean norm of $x \in \Rd$
is denoted by $\|x\|$, and $\|x\|_\infty$ denotes
its supremum norm.
For $f:\R^d \rightarrow \R$ let
\[
\|f\|_\infty = \sup_{x \in \R^d} |f(x)|
\]
denote its supremum norm.
A finite collection $f_1, \dots, f_N:\Rd \rightarrow \R$
  is called an $\varepsilon$--$L_1$-- cover of $\F$
on $x_1^n =(x_1,\dots,x_n)  \in (\Rd)^n$
  if for any $f \in \F$ there exists  $i \in \{1, \dots, N\}$
  such that
  \[
  \frac{1}{n}
\sum_{j=1}^n
|f(x_j)-f_i(x_j)|
< \varepsilon.
  \]
  The $\varepsilon$--$L_1$- covering number of $\F$ on $x_1^n$
  is the  size $N$ of the smallest $\varepsilon$--$L_1$-- cover
  of $\F$ on $x_1^n$
and is denoted by $\Nu_1(\varepsilon,\F, x_1^n)$.

\subsection{Outline}
  The outline of this paper is as follows: In Section \ref{se2} we
  define
our neural network regression estimates and in Section \ref{se3} we
present
our main theoretical result.
The finite sample size performance of our newly proposed estimate
is illustrated in Section \ref{se4} by applying it to simulated
data. The proofs
are given in Section \ref{se5}.

\section{Definition of the estimate}
\label{se2}
In the construction of our estimate we assume that
the regression function $m$ satisfies (\ref{se1eq2})
and that the support of $X$ is contained in the cube
$[-A,A]^d$ for some given $A \geq 1$.
We approximate each $g_s:\R \rightarrow \R$ by a neural
network with logistic squasher
\[
\sigma(x)=\frac{1}{1+e^{-x}}
\]
chosen such that it is close to a piecewise
constant function of the form
\[
u \mapsto \sum_{l=1}^K a_{s,l} \cdot 1_{[b_l, \infty)} + a_{s,0}.
\]
As we will show in Lemma \ref{le5} below, such a neural network
can be chosen of the form
\[
u \mapsto \sum_{l=1}^K a_{s,l} \cdot \sigma( \rho_n \cdot (u-b_l)) + a_{s,0},
\]
where $\rho_n>0$ is a large constant,
and the error of this approximation will be small at all those points,
where $\rho_n \cdot |u-b_l|$ is large. By replacing $u$ with $\bc_s^T x$
we see that we can approximate $m$ by networks
with one hidden layer and $K \cdot r$ neurons in this hidden layer
defined by
\begin{equation}
  \label{se2eq1}
f_{net,(\ba,\bb)}(x)=\sum_{k=1}^{K \cdot r} a_k \cdot \sigma
\left(
\sum_{j=1}^d b_{k,j} \cdot x^{(j)} + b_{k,0}
\right)+ a_0.
\end{equation}
Here, $K \cdot r\in \N$ is the number of neurons, $\sigma:\R \rightarrow \R$
is the activation function and 
\[
a_k \in \R \quad (k=0,\dots,K \cdot r) \quad \mbox{and} \quad
b_{k,j} \in \R \quad (k=1,\dots,K \cdot r, j=0, \dots, d)
\]
are the weights. The above condition that $\rho_n \cdot |u-b_l|$
is large in order to achieve a small error at point $u$ of the above
neural network approximation of the piecewise constant function
is replaced by the assumption that
\[
\min_{i=1,\dots,n}
| \sum_{j=1}^d b_{k,j} \cdot X_i^{(j)} + b_{k,0}|
\]
is large, which will enable us to show that our approximation is
good at all $x$-values of the data points. And this condition in turn
will be ensured by a proper choice of the initial weights described below.

We will learn the weights by gradient descent. More precisely, we minimize
the penalized empirical $L_2$ risk
\begin{equation}
\label{se2eq2}
F(\ba,\bb)=
\frac{1}{n} \sum_{i=1}^n
| f_{net,(\ba,\bb)}(X_i)-Y_i|^2
+
\frac{c_1}{n}
\cdot \sum_{k=0}^{K \cdot r} a_k^2,
\end{equation}
where $c_1 >0$ is a constant,
by choosing an appropriate starting value $(\ba^{(0)},\bb^{(0)})$ and
by setting
\begin{equation}
\label{se2eq3}
\left(
{\ba^{(t+1)} \atop \bb^{(t+1)} }
\right)
=
\left(
{\ba^{(t)} \atop \bb^{(t)} }
\right)
- \lambda_n
\cdot
(\nabla_{(\ba,\bb)} F)(\ba^{(t)},\bb^{(t)})
\end{equation}
for some $\lambda_n>0$ chosen below and $t=0,1,\dots,t_n-1$.

Next, we explain how we choose the initial values
$(\ba^{(0)},\bb^{(0)})$
for our weights.
As explained above,
our choice is motivated by the structure of $m$ in
 the
projection pursuit model (\ref{se1eq2}). Here the number
$r$ of terms in this model is a parameter of our estimate (which we
will choose data-dependent in any application, cf.,  Remark 2 below).
In a first step we randomly choose values
\begin{equation}
\label{se2eq4}
\bar{\bc}_1, \dots, \bar{\bc}_r \in [-1,1]^d
\end{equation}
as an independent sample from a uniform distribution on $[-1,1]^d$ such that $\|\bar{\bc}_s\| = 1$ $(s = 1, \dots, r)$.
Using these values as approximation of the directions
$\bc_1, \dots, \bc_r$ of our projection pursuit model,
we define our initial inner weights as follows: 
For $s \in \{1,\dots,r\}$ we define
\[
b_{(s-1) \cdot K + 1,0}, \dots, b_{(s-1) \cdot K + 1,d}, 
\dots, b_{s\cdot K,0 }, \dots, b_{s\cdot K,d } 
\]
according to $\bar{\bc}_s$ and to $X_1,\dots,X_n$: First, we choose $b_1, \dots, b_K \in
\R$ such that $b_1 < b_2 < \dots < b_K$ and
\[
b_1 \leq -A \cdot \sqrt{d},
\]
\[
b_K \geq A \cdot \sqrt{d} - \frac{4 \cdot \sqrt{d} \cdot A}{K-1},
\]
\[
\frac{\sqrt{d} \cdot A}{(n+1) \cdot (K-1)}
\leq
|b_{k+1}-b_k| 
\leq 
\frac{4 \cdot \sqrt{d} \cdot A}{K-1} 
\quad (k=1,\dots,K-1)
\]
and
\[
\min_{i=1, \dots, n, k=1, \dots, K}
\left|
\bar{\bc}^T_s X_i- b_k
\right|
\geq
\frac{\sqrt{d} \cdot A}{(n+1) \cdot (K-1)}.
\]
Such a choice is always possible, e.g., we can
set $b_1 = -\sqrt{d} \cdot A - 2 \cdot \sqrt{d} \cdot A/((n+1) \cdot (K-1))$
and define $b_k$ $(k=2, \dots, K)$
by subdividing the interval 
\[
\left[
-\sqrt{d} \cdot A + (k-2) \cdot \frac{2 \cdot \sqrt{d} \cdot A}{K-1}, 
-\sqrt{d} \cdot A + (k-1) \cdot \frac{2 \cdot \sqrt{d} \cdot A}{K-1}
\right]
\]
into $(n+1)$ equidistant subintervals of length $2 \cdot \sqrt{d} \cdot A/((K-1) \cdot (n+1))$
and by choosing $b_k$ as the midpoint of one of those
intervals which does not contain any of the $n$ values
$\bar{\bc}^T_s X_i$ (such an interval must exist since not every one of the
$n+1$ disjoint intervals can contain one of the above $n$ points).
As soon as we have chosen $b_1, \dots, b_K$ we define $b_{(s-1) \cdot K +k,j}$
($s=1, \dots, r$, $k=1, \dots, K$, $j=0, \dots, d$)
such that
we have for some $\rho_n>0$ chosen below (cf., Theorem \ref{th1} below)
\[
\sum_{j=1}^d b_{(s-1) \cdot K +k,j} \cdot x^{(j)} + b_{(s-1) \cdot K +k,0}
=
\rho_n \cdot (\bar{\bc}_s^T x - b_k)
\quad \mbox{for all } x \in \Rd,
\]
namely, we set
\[
b_{(s-1) \cdot K +k,j} = \rho_n \cdot \bar{\bc}_s^{(j)} \quad \mbox{and} \quad
b_{(s-1) \cdot K +k,0}=-\rho_n \cdot b_k
\]
($s=1, \dots, r$, $k=1, \dots, K$, $j=1, \dots, d$).
Then, we choose $a_l=0$ for all $l \in \{0, \dots, K \cdot r\}$.

After this choice of $(\ba^{(0)},\bb^{(0)})$ we define
$(\ba^{(t+1)},\bb^{(t+1)})$ 
recursively by (\ref{se2eq3})
for $\lambda_n>0$ and $t=0,1,\dots,t_n-1$.

We repeat this whole procedure $I_n$ times, and let
\[
\tilde{m}_n
\]
be the neural network which achieves the smallest penalized
empirical $L_2$ error (\ref{se2eq2}) among all the $I_n$ networks.
Finally we truncate our estimate by selecting some $\beta_n>0$ and
by setting
\[
m_n(x) = T_{\beta_n} \tilde{m}_n(x),
\]
where
$T_{\beta_n} z = \max\{ \min\{ z, \beta_n\},-\beta_n\}$
for $z \in \R$.

\section{Main result}
\label{se3}
Our main result is the following theorem.
\begin{theorem}
\label{th1}
Let $n \geq 2$, let $A \geq 1$ and let $(X,Y)$, $(X_1,Y_1)$, \dots,
$(X_n,Y_n)$ be independent and identically distributed
random variables with values in $[-A,A]^d \times \R$.
Set $m(x)=\EXP\{Y | X=x\}$ and assume that $(X,Y)$ satisfies
\begin{align}\label{subgaus}
\mathbf E \left( e^{c_{2}\cdot |Y |^2}\right) <\infty
\end{align}
for some constant $c_{2}>0$, and that
$m$ satisfies
\[
m(x)=\sum_{s=1}^r g_s(\bc_s^T x) \quad (x \in \Rd)
\]
for some $r \in \N$, $\bc_s \in [-1,1]^d$, where $\|\bc_s\| = 1$, and $g_s:\R \rightarrow \R$ $(s=1,
\dots,r)$. Assume that $g_s$ is $(p,C)$-smooth for $s \in \{1, \dots, r\}$,
where $p \in (0,1]$ and $C>0$ are fixed.
Define the regression estimate $m_n$ as in Section \ref{se2}
with
\[
\sigma(x)=\frac{1}{1+e^{-x}},
\]
with parameter $r$ as in the above projection pursuit model, and with
the other parameters chosen by
\[
\beta_n=c_3 \cdot \log n, \quad
K=K_n= \lceil (n/ (\log n)^3)^{1/(2p+1)} ]\rceil, \quad
  \lambda_n = \frac{1}{3\cdot K \cdot r }, \quad
\rho_n=n^2 \cdot K, 
\]
and
\[
t_n=K_n \cdot n \cdot (\log n)^2
\quad \mbox{and} \quad
I_n = \lceil (\log n)^{-3 \cdot r \cdot d/(2p+1)} \cdot n^{r \cdot d/(2p+1)} \rceil. 
\]
Then $m_n$ satisfies
\[
\EXP \int |m_n(x)-m(x)|^2 \PROB_X (dx)
\leq
c_4 \cdot
\left(
\frac{(\log n)^3}{n}
\right)^{\frac{2p}{2p+1}}
\]
for some constant $c_4>0$ which does not depend on $n$. 
\end{theorem}

\noindent
{\bf Remark 1.}
According to Stone (1982) the rate of convergence in  the above
theorem is optimal up to a logarithmic factor in case
of a $(p,C)$-smooth projection pursuit model. Because of the
fact that this rate of convergence is independent of the dimension $d$
of $X$, the above theorem shows that our newly proposed
computable neural network regression estimate is able
to circumvent the curse of dimensionality in case that the regression
function satisfies the assumption of projection pursuit.
We should however mention that the number of repitions $I_n$
of the initial  random choices of the directions $\bar{\bc}_s$
and correspondingly the number of repititions of the $t_n$ gradient
descent steps is rather huge.

\noindent
{\bf Remark 2.}
The parameters $r$ and $K_n$, and also $I_n$, of the above algorithm
depend on the projection pursuit model and hence are unknown
in any application. However, it is easy to choose them data-dependently
by using, e.g., the splitting of the sample technique as explained in
the next section. In this way it is possible to define an estimate
which
does not depend on the value of $r$ of the projection pursuit model
and which is nevertheless able to achieve the rate of convergence
in Theorem \ref{th1}.

\section{Application to simulated data}
\label{se4}
In this section we illustrate the finite sample size performance of
our newly proposed
estimate by applying it
to simulated data.

The simulated data which we use is defined as follows: We choose
$d=4$, $X$ uniformly distributed on $[-1,1]^d$, $\epsilon$ standard
normal and independent of $X$, and we define $Y$ by
\begin{equation}
  \label{se4eq1}
  Y=m_j(X) + \sigma \cdot \tau_j \cdot \epsilon,
  \end{equation}
where $m_j:[-1,1]^d \rightarrow \R$ is described below,
$\tau_j>0$ is a scaling value defined below and $\sigma$
is chosen from $\{0.05,0.2\}$ $(j \in \{1,2\})$.
As regression function we use 
\[
m_1(x_1,x_2,x_3,x_4)=2 \cdot \sin \left(
\frac{2 \cdot \pi}{\sqrt{4}} \cdot (-x_1+x_2-x_3+x_4)
\right),
\]
so $m_1$ satisfies a single index model, and
\begin{eqnarray*}
  &&
  m_2(x_1,x_2,x_3,x_4)\\
  &&
  =
  4 \cdot \sin \left(
\frac{2 \cdot \pi}{\sqrt{4}} \cdot (-x_1+x_2-x_3+x_4) 
\right)
+
\frac{7}
     {2+ \frac{1}{\sqrt{30}} \cdot (x_1 - 2 \cdot x_2 + 3 \cdot x_3 - 4 \cdot x_4)},
\end{eqnarray*}
hence $m_2$ satisfies a single index model with $r=2$ terms.
$\tau_j$ is chosen approximately as IQR of samples of size
$100,000$ of $m(X)$,
and we use the concrete values $\tau_1=2.8289$ and $\tau_2=5.2841$.
From this distribution we generate samples of size $n=100$ and $n = 200$ and apply our
newly proposed neural network regression estimate and two alternative
regression estimates to these samples. Then we compute the $L_2$ errors
of these three estimates approximately by using the empirical $L_2$
error
$\varepsilon_{L_2,\bar{N}}(\cdot)$
on an independent sample of $X$ of size $\bar{N}=10,000$.
Since this error strongly depends on the behavior of the correct function $m_j$, we consider it in relation to the error of the simplest estimate for $m_j$ we can think of, a completely constant function (whose value is the average of the observed data according to the least squares approach). Thus, the scaled error measure we use for evaluation of the estimates is $\varepsilon_{L_2,\bar{N}}(m_{n,i})/\bar{\varepsilon}_{L_2,\bar{N}}(avg)$, where $\bar{\varepsilon}_{L_2,\bar{N}}(avg)$ is the median of $50$ independent realizations of the value
one obtains if one plugs
the average of $n$ observations into $\varepsilon_{L_2,\bar{N}}(\cdot)$. To a certain extent, this quotient can be interpreted as the relative part of the error of the constant estimate that is still contained in the more sophisticated approaches.
The resulting scaled
errors of course depend on the random sample of $(X,Y)$, and to
be able to compare these values nevertheless we repeat the whole
computation $25$ times and report the median and the interquartile
range of the $25$ scaled errors for each of our three estimates.

Our first estimate {\it Tps} is a smoothing spline estimate with parameter
chosen by generalized cross validation as implemented in the
routine {\it Tps()} of the library {\it fields} in {\it R}.

Our second estimate {\it neighbor}
is a nearest neighbor estimate where the number
of nearest neighbors is chosen from the set $\{1,2,4,8,16,32\}$
by splitting of the sample. Here we split our sample in a learning
sample of size $n_l=0.8 \cdot n$ and a testing sample of size
$n_t=0.2 \cdot n$. We compute the estimate for all parameter
values from the above set using the learning sample, compute
the corresponding empirical $L_2$ risk on the testing sample
and choose the parameter value which leads to the minimal
empirical $L_2$ risk on the testing sample.

Our third estimate {\it neural} is our newly proposed
neural network estimate presented in this paper, which
we have implemented in {\it R}.
Here the parameters $r$ and $K$ of the estimate are chosen
via splitting of the sample
(as described above) from the set $\{1,2\}$ and $\{5,10,20\}$,
respectively. In order
to accelerate the computation of this estimate we use
only $I_n=50$ random choices for the vectors of directions
in the computation of the estimate for each parameter
value.

The results are summarized in Table \ref{se4ta1}
and in Table \ref{se4ta2}.
As we can see from the reported scaled errors, our
newly proposed neural network estimate outperforms
in both cases
in all four settings both the smoothing spline estimate
and the nearest neighbor estimate.

\begin{table}[htbp]                                  
\centering
\makebox[\textwidth][c]{                                     
\begin{tabular}{|c|c|c|c|c|}
\hline
&\multicolumn{2}{|c|}{$m_1$} & \multicolumn{2}{|c|}{$m_2$}\\
\hline
\textit{noise} & {$5\%$} & {$20\%$} & {$5\%$} & {$20\%$}\\
\hline
$\bar{\varepsilon}_{L_2,\bar{N}}(avg)$ & {$2.0154$} &  {$2.0219$} &  {$10.3521$} &  {$10.3627$} \\
\hline                                         
\textit{approach} & median (IQR) & median (IQR) & median (IQR) & median (IQR) \\                  
\hline
Tps & $1.18$ $(0.17)$  & $1.19$ $(0.14)$  & $0.89$ $(0.09)$ & $0.98$  $(0.17)$ \\
neighbor & $1.06$ $(0.15)$ & $1.13$ $(0.27)$   & $0.91$ $(0.07)$  & $0.92$  $(0.08)$ \\
neural & $0.52$ $(0.34)$  & $0.46$ $(0.25)$   & $0.42$ $(0.16)$  & $0.56$ $(0.15)$   \\
\hline
\end{tabular}
}
\caption{Median and IQR of the scaled empirical $L_2$ error of estimates for $m_1$ and $m_2$ for sample size $n=100$.}                
\label{se4ta1}                            
\end{table}

\begin{table}[htbp]                                  
	\centering
	\makebox[\textwidth][c]{                                     
		\begin{tabular}{|c|c|c|c|c|}
			\hline
			&\multicolumn{2}{|c|}{$m_1$} & \multicolumn{2}{|c|}{$m_2$}\\
			\hline
			\textit{noise} & {$5\%$} & {$20\%$} & {$5\%$} & {$20\%$}\\
			\hline
			$\bar{\varepsilon}_{L_2,\bar{N}}(avg)$ & {$2.0125$} &  {$2.0109$} &  {$10.3127$} &  {$10.3192$} \\
			\hline                                         
			\textit{approach} & median (IQR) & median (IQR) & median (IQR) & median (IQR) \\                  
			\hline
			Tps & $0.75$ $(0.08)$  & $0.82$ $(0.17)$  & $0.55$ $(0.05)$ & $0.64$  $(0.08)$ \\
			neighbor & $0.88$ $(0.08)$ & $0.96$ $(0.08)$   & $0.70$ $(0.07)$  & $0.77$  $(0.10)$ \\
			neural & $0.44$ $(0.29)$  & $0.44$ $(0.31)$   & $0.34$ $(0.19)$  & $0.40$ $(0.22)$   \\
			\hline
		\end{tabular}
	}
	\caption{Median and IQR of the scaled empirical $L_2$ error of estimates for $m_1$ and $m_2$ for sample size $n=200$.}                
	\label{se4ta2}                            
\end{table}

\section{Proofs}
\label{se5}

\subsection{Learning of linear penalized least squares estimates by gradient descent}

Let $(x_1,y_1), \dots, (x_n,y_n) \in \Rd \times \R$, let $K \in \N$,
let $B_1,\dots,B_K:\Rd \rightarrow \R$ and let $c_1>0$. In this subsection
we consider the problem to minimize
\begin{equation}
  \label{se5eq1}
  F(\ba) =
  \frac{1}{n} \sum_{i=1}^n
  | \sum_{k=1}^K a_k \cdot B_k(x_i)-y_i|^2
  +
  \frac{c_1}{n} \cdot \|\ba\|^2,
  \end{equation}
where
\[
\ba=(a_1,\dots,a_K)^T \quad \mbox{and} \quad
\|\ba\|^2=\sum_{j=1}^K a_j^2,
\]
by gradient descent. To do this, we choose $\ba^{(0)} \in \R^K$
and set
\begin{equation}
  \label{se5eq2}
  \ba^{(t+1)} = \ba^{(t)}
  - \lambda_n \cdot (\nabla_\ba F)(\ba^{(t)})
    \end{equation}
for some properly chosen $\lambda_n>0.$

\begin{lemma}
  \label{le1}
  Let $F:\R^K \rightarrow \R$ be a differentiable function and define
  $\ba^{(t+1)}$ by (\ref{se5eq2}),
  where
  \begin{equation}
    \label{le1eq1}
    \lambda_n=\frac{1}{L_n}
  \end{equation}
  for some $L_n>0$. Let $\ba_{opt} \in \R^K$ be arbitrary.

  \noindent
      {\bf a)} If
      \begin{equation}
        \label{le1eq2}
        \| (\nabla_\ba F)(\ba_1) - (\nabla_\ba F)(\ba_2)\|
        \leq L_n \cdot \|\ba_1-\ba_2\| \quad (\ba_1,\ba_2 \in \R^K)
      \end{equation}
      holds, then we have
      \[
      F(\ba^{(t+1)})-F(\ba^{(t)}) \leq - \frac{1}{2 \cdot L_n}
      \cdot \|(\nabla_\ba F)(\ba^{(t)})\|^2.
      \]

        \noindent
            {\bf b)} If inequality (\ref{le1eq2}) and, in addition,
      \begin{equation}
        \label{le1eq3}
        \|(\nabla_\ba F)(\ba)\|^2
        \geq \rho_n \cdot (F(\ba)-F(\ba_{opt})) \quad (\ba \in \R^K)
      \end{equation}
      hold for some $\rho_n > 0$, then we have
      \[
      F(\ba^{(t+1)})-F(\ba_{opt})
      \leq
      \left(1 - \frac{\rho_n}{2 \cdot L_n} \right)
      \cdot
(      F(\ba^{(t)})-F(\ba_{opt})
).      
      \]
  \end{lemma}

\noindent
    {\bf Proof.} Lemma \ref{le1} follows
    from well-known bounds in the literature, see, e.g.,
    Karimi, Nutini and Schmidt (2018). For the sake
    of completeness a complete proof is given in the
    supplementary material. \hfill $\Box$

        \begin{lemma}
          \label{le2}
          Let $F$ be defined by (\ref{se5eq1}). Then we have for any
          $\ba_1,\ba_2 \in \R^K$
          \begin{eqnarray*}
            &&
        \| (\nabla_\ba F)(\ba_1) - (\nabla_\ba F)(\ba_2)\|
        \leq
        \left(
        2 \cdot \sum_{k=1}^K \frac{1}{n} \sum_{i=1}^n B_k(x_i)^2
        +
        \frac{2 \cdot c_1}{n}
        \right)
        \cdot \|\ba_1-\ba_2\|.            
            \end{eqnarray*}
        \end{lemma}

        \noindent
            {\bf Proof.}
            We have
            \[
F(\ba) = \frac{1}{n} \cdot \left( \bB \cdot \ba - \by \right)^T \cdot \left( \bB \cdot \ba - \by \right) + \frac{c_1}{n} \cdot \ba^T \cdot \ba
\]
where
\[
\bB = \left(
B_j(x_i)
\right)_{1 \leq i \leq n, 1 \leq j \leq K}
\quad \mbox{and} \quad
\by=(y_1, \dots, y_n)^T.
\]
Consequently,
\[
(\nabla_\ba F)(\ba) = \frac{2}{n} \cdot \left( \bB^T \bB \ba - \bB^T \by \right)
+
\frac{2 \cdot c_1}{n} \cdot \ba
\]
and
\[
      \| (\nabla_\ba F)(\ba_1) - (\nabla_\ba F)(\ba_2)\|
        \leq
        \|  \frac{2}{n} \cdot \bB^T \bB \cdot (\ba_1-\ba_2) \|
        +
        \frac{2 \cdot c_1}{n} \cdot \| \ba_1-\ba_2 \|.
        \]
        By applying twice the Cauchy-Schwarz inequality we get
        \begin{eqnarray*}
          \left\|  \frac{2}{n} \cdot \bB^T \bB \cdot \ba \right\|^2
          &
          =&
          \sum_{j=1}^K \left(
\sum_{k=1}^K (\frac{2}{n} \sum_{i=1}^n B_j(x_i) \cdot B_k(x_i)) \cdot a_k
\right)^2
\\
&
\leq &
          \sum_{j=1}^K 
\sum_{k=1}^K (\frac{2}{n} \sum_{i=1}^n B_j(x_i) \cdot B_k(x_i))^2 \cdot 
\|\ba\|^2
\\
&
\leq &
\sum_{j=1}^K 
\sum_{k=1}^K 4 \cdot \frac{1}{n} \sum_{i=1}^n B_j(x_i)^2 \cdot
\frac{1}{n} \sum_{i=1}^n
B_k(x_i))^2 \cdot 
\|\ba\|^2
\\
&
= &
\left(
2 \cdot
\sum_{k=1}^K \frac{1}{n} \sum_{i=1}^n B_k(x_i)^2
\right)^2
\cdot 
\|\ba\|^2,
          \end{eqnarray*}
        which implies the assertion.
            \hfill $\Box$

        \begin{lemma}
          \label{le3}
          Let $F$ be defined by (\ref{se5eq1}) and choose $\ba_{opt}$
          such that
          \[
F(\ba_{opt})=\min_{\ba \in \R^K} F(\ba).
          \]
          Then for any
          $\ba \in \R^K$
          we have 
                    \[
                            \|(\nabla_\ba F)(\ba)\|^2
                            \geq \frac{4 \cdot c_1}{n} \cdot (F(\ba)-F(\ba_{opt})).
                            \]
          \end{lemma}

        \noindent
            {\bf Proof.} Set
            \[
\bB = \left(
B_j(x_i)
\right)_{1 \leq i \leq n, 1 \leq j \leq K}
\quad \mbox{and} \quad
\bA=\frac{1}{n} \cdot \bB^T \cdot \bB + \frac{c_1}{n} \cdot \bID,
            \]
            where $\bID$ is the unit matrix. Then $\bA$ is positive definite
            and hence regular, from which we can conlcude
            \begin{eqnarray*}
              F(\ba)         & = &
              \frac{1}{n} \cdot \left( \bB \cdot \ba - \by \right)^T \cdot \left( \bB \cdot \ba - \by \right) + \frac{c_1}{n} \cdot \ba^T \cdot \ba
              \\
              &=&
              \ba^T \bA \ba - 2 \by^T \frac{1}{n} \bB \ba + \frac{1}{n} \by^T \by
              \\
              &=&
              (\ba - \bA^{-1} \frac{1}{n} \bB^T \by)^T \bA (\ba - \bA^{-1} \frac{1}{n} \bB^T \by)
              + F(\ba_{opt}),
              \end{eqnarray*}
            where
            \[
            F(\ba_{opt})
            =
            \frac{1}{n} \by^T \by
            -
            \by^T \cdot \frac{1}{n} \cdot \bB \bA^{-1} \cdot \frac{1}{n} \cdot
            \bB^T \by.
            \]
            Using
            \[
\bb^T \bA \bb \geq \frac{c_1}{n} \cdot \bb^T \bb
            \]
   and $\bA^T=\bA$         we conclude
   \begin{eqnarray*}
     &&
     F(\ba)-F(\ba_{opt})
     \\
     &&
     =
     ((\bA^{1/2})^T  (\ba - \bA^{-1} \frac{1}{n} \bB^T \by))^T \bA^{1/2} (\ba - \bA^{-1} \frac{1}{n} \bB^T \by)
     \\
     &&
     \leq
     \frac{n}{c_1}
     \cdot
     ((\bA^{1/2})^T  (\ba - \bA^{-1} \frac{1}{n} \bB^T \by))^T \bA \bA^{1/2} (\ba - \bA^{-1} \frac{1}{n} \bB^T \by)
     \\
     &&
     =
    \frac{n}{c_1}
     \cdot
          ((\bA)^T  (\ba - \bA^{-1} \frac{1}{n} \bB^T \by))^T \bA  (\ba - \bA^{-1} \frac{1}{n} \bB^T \by)
  \\
     &&
     =
    \frac{n}{c_1}
     \cdot
          (\bA  \ba - \frac{1}{n} \bB^T \by)^T (\bA \ba -  \frac{1}{n} \bB^T \by)\\
     &&
     =
    \frac{n}{4 \cdot c_1}
     \cdot
          (2 \bA  \ba - \frac{2}{n} \bB^T \by)^T (2 \bA \ba -  \frac{2}{n} \bB^T \by)\\
     &&
     =
    \frac{n}{4 \cdot c_1}
     \cdot
   \left\| (\nabla_\ba F)(\ba) \right\|^2 , 
   \end{eqnarray*}
   where the last equality follows from
   \begin{eqnarray*}
     &&
     (\nabla_\ba F)(\ba)
     =
     \nabla_\ba \left(
              \ba^T \bA \ba - 2 \by^T \frac{1}{n} \bB \ba + \frac{1}{n} \by^T \by
              \right)
              =
2 \bA  \ba - \frac{2}{n} \bB^T \by.              
     \end{eqnarray*}
            \quad
            \hfill $\Box$

\subsection{Result for neural networks with one hidden layer}
In this subsection we study neural networks with one
hidden layer, which are defined by
\begin{equation}
  \label{se5sub2eq1}
f_{net,(\ba,\bb)}(x)=\sum_{k=1}^K a_k \cdot \sigma
\left(
\sum_{j=1}^d b_{k,j} \cdot x^{(j)} + b_{k,0}
\right)+ a_0
\end{equation}
(compare (\ref{se2eq1})),
where $K \in \N$ is the number of neurons, $\sigma:\R \rightarrow \R$
is the activation function and where the weights
\[
a_k \quad (k=0,\dots,K) \quad \mbox{and} \quad
b_{k,j} \in \R \quad (k=1,\dots,K, j=0, \dots, d)
\]
are learned by gradient descent. More precisely, we minimize
\begin{equation}
\label{se5sub2eq2}
F(\ba,\bb)=
\frac{1}{n} \sum_{i=1}^n
| f_{net,(\ba,\bb)}(x_i)-y_i|^2
+
\frac{c_1}{n}
\cdot \sum_{k=0}^K a_k^2
\end{equation}
(compare (\ref{se2eq2}))
by choosing an appropriate starting value $(\ba^{(0)},\bb^{(0)})$ and
by setting
\begin{equation}
\label{se5sub2eq3}
\left(
{\ba^{(t+1)} \atop \bb^{(t+1)} }
\right)
=
\left(
{\ba^{(t)} \atop \bb^{(t)} }
\right)
- \lambda_n
\cdot
(\nabla_{(\ba,\bb)} F)(\ba^{(t)},\bb^{(t)})
\end{equation}
for some $\lambda_n>0$ chosen below.

Our main idea is, that in the case of the logistic squasher
\[
\sigma(x)=\frac{1}{1+e^{-x}} \quad (x \in \R),
\]
the neural network (\ref{se5sub2eq1})
is
for appropriate weigths
$b_{k,j}$
close to a linear combination of indicator functions,
and in this case the gradient descent will change
the inner weights
$b_{k,j}$ only slightly. From this we will conclude
from our results for linear least squares estimates that
for such networks the gradient descent leads to estimates
where the outer weights $a_k$ are chosen optimally.

In Lemma \ref{le5} below we study the approximation of
H\"older continuous functions by neural networks of the above
form in the case of univariate functions and networks. To do this,
we will need the following auxiliary result.

\begin{lemma}
  \label{le4}
  Let $\sigma$ be the logistic squasher.

    \noindent
  {\bf a)} For any $x \in \R$ we have
  \[
|\sigma(x)-1_{[0,\infty)}(x)| \leq e^{-|x|}.
  \]

    \noindent
  {\bf b)} For any $b \in \R$, $c>0$ and $x \in \R$ we have
  \[
  | \sigma( c \cdot (x-b)) - 1_{[b,\infty)}(x)|
    \leq
    e^{-c \cdot |x-b|}.
  \]
  \end{lemma}

\noindent
    {\bf Proof.} {\bf a)} For $x \geq 0$ we have
    \[
    |\sigma(x)-1_{[0,\infty)}(x)| = 1 - \frac{1}{1+e^{-x}} = \frac{e^{-x}}{1+e^{-x}}
      \leq e^{-x}=e^{-|x|}.
    \]
    And for $x<0$ we get
    \[
    |\sigma(x)-1_{[0,\infty)}(x)| = \frac{1}{1+e^{-x}}
      \leq e^{x}=e^{-|x|}.
    \]

    \noindent
    {\bf b)}
    From $c>0$ and  a)  we get
    \[
    | \sigma( c \cdot (x-b)) - 1_{[b,\infty)}(x)|
      =
    | \sigma( c \cdot (x-b)) - 1_{[0,\infty)}(c \cdot (x-b))|      
    \leq
    e^{-|c \cdot (x-b)|}
    =
    e^{-c \cdot |x-b|}.
    \]
    \hfill $\Box$

\begin{lemma}
  \label{le5}
  Let $\sigma$ be the logistic squasher.
  Let $\bar{\bc} \in [-1,1]^d$ with $\| \bar{\bc} \| = 1$ and let $g:\R \rightarrow \R$
  be $(p,C)$-smooth for some $p \in (0,1]$ and $C>0$.
  Let $\rho_n > 0$, $K \in \N$ 
  and choose $b_1,b_2, \dots, b_K \in \R$ such that $b_1 < b_2 < \dots < b_K$ and
\[
b_1 \leq -A \cdot \sqrt{d},
\]
\[
b_K \geq A \cdot \sqrt{d} - \frac{4 \cdot A \cdot \sqrt{d}}{K-1}
\]
and
\[
\frac{A \cdot \sqrt{d}}{(n+1)\cdot(K-1)}
\leq
|b_{k+1}-b_k| 
\leq 
\frac{4 \cdot A \cdot \sqrt{d}}{K-1} 
\quad (k=1,\dots,K-1).
\]
Let
\[
a_0 = g(b_1) \quad \mbox{and} \quad
a_k = g(b_k)-g(b_{k-1}) \quad (k=1, \dots, K ).
\]
Then we have
\begin{eqnarray*}
&&
\sup_{x \in [-A,A]^d}
\left| a_0 + \sum_{k=1}^{K} a_k \cdot \sigma(\rho_n \cdot (\bar{\bc}^T x - b_k)) - g(\bar{\bc}^T x) \right|
\\
&&
\leq
 \frac{3 \cdot (4 \cdot A \cdot \sqrt{d})^p \cdot C}{(K-1)^p} + C \cdot (4 \cdot A \cdot \sqrt{d})^p \cdot (K-1)^{1-p}
\cdot e^{- \frac{\rho_n \cdot (A \cdot \sqrt{d})}{(n+1) \cdot (K-1)}}.
\end{eqnarray*}
\end{lemma}

\noindent
    {\bf Proof.}
    We have
    \begin{eqnarray*}
      &&
\left| a_0 + \sum_{k=1}^{K} a_k \cdot \sigma(\rho_n \cdot (\bar{\bc}^T x - b_k)) - g(\bar{\bc}^T x)
\right|
\\
&&
\leq
\left| \sum_{k=1}^{K} a_k \cdot \sigma(\rho_n \cdot (\bar{\bc}^T x - b_k))  -
\sum_{k=1}^{K} a_k \cdot 1_{[b_k,\infty)}(\bar{\bc}^T x)
\right|
\\
&&
\quad
+
\left| a_0 + \sum_{k=1}^{K} a_k \cdot 1_{[b_k,\infty)}(\bar{\bc}^T x) - g(\bar{\bc}^T x)
\right|.
      \end{eqnarray*}
    For $b_j \leq \bar{\bc}^T x < b_{j+1}$, where $j \in \{1,\dots,K-1\}$,
we can conclude from the definition of $a_k$,
from the $(p,C)$-smoothness of $g$
and from our choice of the $b_k$
    \begin{eqnarray*}
      &&
\left| a_0 + \sum_{k=1}^{K} a_k \cdot 1_{[b_k,\infty)}(\bar{\bc}^T x) - g(\bar{\bc}^T x)
  \right|
  \\
  &&
  =
  \left| a_0 + \sum_{k=1}^{j} a_k - g(\bar{\bc}^T x)
  \right|
  =|g(b_{j})-g(\bar{\bc}^T x)| 
  \\
  &&
  \leq C \cdot |b_{j}-\bar{\bc}^T x|^p 
  \leq C \cdot |b_{j+1} - b_j|^p
  \leq \frac{C \cdot (4 \cdot A \cdot \sqrt{d})^p}{(K-1)^p}.
      \end{eqnarray*}
    It is easy to see that
    this inequality is also true for $b_{K} \leq \bar{\bc}^T x \leq \sqrt{d}\cdot A$. Hence, we have shown
    \[
    \sup_{x \in [-A, A]^d}
    \left| a_0 + \sum_{k=1}^{K} a_k \cdot 1_{[b_k,\infty)}(\bar{\bc}^T x) - g(\bar{\bc}^T x)
      \right|
      \leq \frac{C \cdot (4 \cdot A \cdot \sqrt{d})^p}{(K-1)^p}. 
    \]
    We finish the proof by showing
    \begin{eqnarray*}
      &&
   \sup_{x \in [-A, A]^d}
 \left| \sum_{k=1}^{K} a_k \cdot \sigma(\rho_n \cdot (\bar{\bc}^T x - b_k)) -
\sum_{k=1}^{K} a_k \cdot 1_{[b_k,\infty)}(\bar{\bc}^T x)
  \right|
  \\
  &&
  \leq
  \ \frac{2 \cdot (4 \cdot A \cdot \sqrt{d})^p \cdot C}{(K-1)^p} + C \cdot (4 \cdot A \cdot \sqrt{d})^p \cdot (K-1)^{1-p}
  \cdot e^{- \frac{\rho_n \cdot (A \cdot \sqrt{d})}{(n+1) \cdot (K-1)}}.
    \end{eqnarray*}
    For $b_j \leq \bar{\bc}^T x \leq b_{j+1}$, where $j \in \{1, \dots, K-1\}$, we have
    \begin{eqnarray*}
      &&
 \left| \sum_{k=1}^{K} a_k \cdot \sigma(\rho_n \cdot (\bar{\bc}^T x-b_k)) -
\sum_{k=1}^{K} a_k \cdot 1_{[b_k,\infty)}(\bar{\bc}^T x)
  \right|
  \\
  &&
  \leq
  \sum_{k=1}^{j-1} |a_k| \cdot \left|  \sigma(\rho_n \cdot (\bar{\bc}^T x-b_k)) -1_{[b_k,\infty)}(\bar{\bc}^T x)
    \right|
    +|a_j| + |a_{j+1}|
    \\
    &&
    \quad
    +
   \sum_{k=j+2}^{K} |a_k| \cdot \left|  \sigma(\rho_n \cdot (\bar{\bc}^T x-b_k)) -1_{[b_k,\infty)}(\bar{\bc}^T x)
      \right|
      \\
      &&
      \leq
      \max_{k=1, \dots, K} |a_k| \cdot
      \left(
      2 + (K-2) \cdot \max_{k \in \{1,2,\dots,j-1,j+2,j+3,\dots, K\}}
      \left|  \sigma(\rho_n \cdot (\bar{\bc}^T x-b_k)) -1_{[b_k,\infty)}(\bar{\bc}^T x)
      \right|
      \right).
      \end{eqnarray*}
  For $b_{K} \leq \bar{\bc}^T x \leq \sqrt{d}\cdot A$ we get
  \begin{eqnarray*}
  	&&
  	\left| \sum_{k=1}^{K} a_k \cdot \sigma(\rho_n \cdot (\bar{\bc}^T x-b_k)) -
  	\sum_{k=1}^{K} a_k \cdot 1_{[b_k,\infty)}(\bar{\bc}^T x)
  	\right|
  	\\
  	&&
  	\leq
  	\max_{k=1, \dots, K} |a_k| \cdot
  	\left(
  	1 + (K-1) \cdot \max_{k \in \{1,2,\dots, K-1\}}
  	\left|  \sigma(\rho_n \cdot (\bar{\bc}^T x-b_k)) -1_{[b_k,\infty)}(\bar{\bc}^T x)
  	\right|
  	\right).
  \end{eqnarray*}
    By definition of $a_k$ and by the $(p,C)$-smoothness of $g$, we have
    \[
|a_k| \leq C \cdot |b_k -b_{k-1}|^p \leq C \cdot \frac{(4 \cdot A \cdot \sqrt{d})^p}{(K-1)^p},
    \]
    which, together with Lemma \ref{le4}, implies
    for $b_j \leq \bar{\bc}^T x \leq b_{j+1}$, where $j \in \{1, \dots, K-1\}$,
 \begin{eqnarray*}
      &&
 \left| \sum_{k=1}^{K} a_k \cdot \sigma(\rho_n \cdot (\bar{\bc}^T x-b_k)) -
\sum_{k=1}^{K} a_k \cdot 1_{[b_k,\infty)}(\bar{\bc}^T x)
  \right|
  \\
  &&
  \leq
  C \cdot \frac{(4 \cdot A \cdot \sqrt{d})^p}{(K-1)^p} \cdot (2 +
  (K-2) \cdot \max_{k \in \{1,2,\dots,j-1,j+2,j+3,\dots, K\}}
  e^{- \rho_n \cdot |\bar{\bc}^T x-b_k|})
  \\
  &&
  \leq
  \frac{2 \cdot (4 \cdot A \cdot \sqrt{d})^p \cdot C}{(K-1)^p} + C \cdot (4 \cdot A \cdot \sqrt{d})^p \cdot (K-1)^{1-p}
  \cdot e^{- \frac{\rho_n \cdot (A \cdot \sqrt{d})}{(n+1) \cdot (K-1)}}.
      \end{eqnarray*}
  It is easy to see that this bound is also true for $b_{K} \leq \bar{\bc}^T x \leq \sqrt{d}\cdot A$. This concludes the proof.
\hfill $\Box$

\begin{lemma}
  \label{le6}
  Let $\sigma$ be the logistic squasher.
Define $F$ by (\ref{se5sub2eq2}) and set
\[
\bar{\bb}=\bb - \lambda_n \cdot (\nabla_\bb F)(\ba,\bb)
\]
for some $\lambda_n>0$, where
\[
\ba = (a_1, \dots, a_K)^T \in \R^{K}
\quad
\text{and}
\quad
\bb = (b_{1,0}, b_{1,1}, \dots, b_{1,d},\dots,b_{K,0},b_{K,1}\dots,b_{K,d})^T \in \R^{K \cdot (d+1)} 
.
\]
Then we have for any $k \in \{1, \dots,K\}$ and any
$j \in \{0,\dots,d\}$:
\begin{eqnarray*}
|\bar{b}_{k,j}-b_{k,j}|
&
\leq &
\lambda_n \cdot 2 \cdot \sqrt{F(\ba,\bb)} \cdot
\max\{1, \max_{i,l} \{|x_{i}^{(l)} | \} \}
\cdot
|a_k|
\\
&&
\hspace*{2cm}
\cdot
\exp \left(
-
\min_{i=1,\dots,n}
\left\{
\left|
\sum_{j=1}^d b_{k,j} \cdot x_i^{(j)} + b_{k,0}
\right|
\right\}
\right).
\end{eqnarray*}
\end{lemma}

  \noindent
      {\bf Proof.}    Using
      \[
      |
      \sigma^\prime
      (
      x)
      |
      =
      |\sigma(x) \cdot (1-\sigma(x))|
      \leq
      \min \left\{
|\sigma(x)|,|1-\sigma(x)|
      \right\}
      \leq
      |\sigma(x)-1_{[0,\infty)}(x)|
      \]
      (where the first inequality holds due to $\sigma(x) \in [0,1]$)
      we can conclude from Lemma \ref{le4} that
      \begin{eqnarray*}
        &&
        \max_{i=1,\dots,n}
      	\left|\sigma^\prime \left( \sum_{j=1}^d b_{k,j} \cdot x_{i}^{(j)} + b_{k,0} \right)\right|
        \\
      	&&\leq
        \max_{i=1,\dots,n}
      	\exp\left( - \left| \sum_{j=1}^d b_{k,j} \cdot x_{i}^{(j)} + b_{k,0} \right| \right)
      	\\
      	&&=
      	\exp \left(
      	-
      	\min_{i=1,\dots,n}
      	\left\{
      	\left|\sum_{j=1}^d b_{k,j} \cdot x_i^{(j)} + b_{k,0} \right|
      	\right\}
      	\right).
      \end{eqnarray*}
As a consequence, we get for $k \in \{1, \dots, K\}$ and $j \in \{1,\dots,d\}$ by the Cauchy-Schwarz inequality
\begin{eqnarray*}
&&
\left|
\frac{\partial F}{\partial b_{k,j}}(\ba,\bb)
\right| \\
&&
= 
\left| 
\frac{2}{n}\sum_{i = 1}^{n}(f_{net,(\ba,\bb)}(x_i)-y_i) \cdot
a_k \cdot \sigma^\prime
\left(
\sum_{j=1}^d b_{k,j} \cdot x_{i}^{(j)} + b_{k,0}
\right) \cdot
x_{i}^{(j)}
\right|
\\
&&
\leq
2 \cdot |a_k| \cdot
\frac{1}{n}\sum_{i = 1}^{n} | f_{net,(\ba,\bb)}(x_i)-y_i | \cdot | x_{i}^{(j)}| 
\cdot \left| \sigma^\prime
\left(
\sum_{j=1}^d b_{k,j} \cdot x_{i}^{(j)} + b_{k,0}
\right) \right|
\\
&&
\leq
2 \cdot \sqrt{
\frac{1}{n} \sum_{i=1}^n
| f_{net,(\ba,\bb)}(x_i)-y_i|^2
\cdot (x_i^{(j)})^2
}
\cdot |a_k| \cdot
\sqrt{
\frac{1}{n} \sum_{i=1}^n
|
\sigma^\prime
(
\sum_{j=1}^d b_{k,j} \cdot x_i^{(j)} + b_{k,0}
)
|^2
}
\\
&&
\leq
2 
\cdot \sqrt{ F(\ba,\bb)}
\cdot  \max_{i,l} \{|x_{i}^{(l)}|\}
\cdot |a_k|
\cdot \sqrt{ \frac{1}{n} \sum_{i=1}^n| \sigma^\prime(\sum_{l=1}^d b_{k,l} \cdot x_l^{(j)} + b_{k,0})| ^2
}
\\
&&
\leq
2
\cdot  \sqrt{ F(\ba,\bb)}
\cdot \max_{i,l} \{|x_{i}^{(l)}|\}
\cdot |a_k|
\cdot \exp \left(
	-
	\min_{i=1,\dots,n}
	\left\{
	\left|
	\sum_{j=1}^d b_{k,j} \cdot x_i^{(j)} + b_{k,0}
	\right|
	\right\}
	\right).
\end{eqnarray*}

Hence, we have shown
\begin{eqnarray*}
	&&
|\bar{b}_{k,j}-b_{k,j}|
\\
&&
=
\lambda_n \cdot \left| \frac{\partial F}{\partial b_{k,j}}(\ba,\bb) \right|
\\
&&
\leq
\lambda_n
\cdot 
2
\cdot  \sqrt{ F(\ba,\bb)}
\cdot \max_{i,l} \{| x_{i}^{(l)} |\}
\cdot |a_k|
\cdot \exp \left(
-
\min_{i=1,\dots,n}
\left\{
\left|
\sum_{j=1}^d b_{k,j} \cdot x_i^{(j)} + b_{k,0}
\right|
\right\}
\right)
\end{eqnarray*}
for any $k \in \{1, \dots, K\}$ and any $j \in \{1,\dots,d\}$ .

In case that $k \in \{1, \dots, K\}$ and $j=0$ we get in a similar fashion
\begin{eqnarray*}
|\bar{b}_{k,0}-b_{k,0}|
&=& 
\lambda_n \cdot \left| \frac{\partial F}{\partial b_{k,0}}(\ba,\bb) \right|
\\
&\leq&
\lambda_n
\cdot 
2
\cdot  \sqrt{ F(\ba,\bb)}
\cdot 1
\cdot |a_k|
\cdot \exp \left(
-
\min_{i=1,\dots,n}
\left\{
\left|
\sum_{j=1}^d b_{k,j} \cdot x_i^{(j)} + b_{k,0}
\right|
\right\}
\right),
\end{eqnarray*}
which implies the assertion.
\hfill $\Box$

      \begin{lemma}
  \label{le7}
Define $F$ by (\ref{se5sub2eq2})
 and define $(\ba^{(t)},\bb^{(t)})$ by
(\ref{se5sub2eq3}). Assume that $(\ba^{(t)},\bb^{(t)})$ satisfy
for $t \in \{0,\dots,t_n-1\}$
\begin{equation}
\label{le7eq1}
F(\ba^{(t)},\bb^{(t)}) \leq c_5 < \infty,
\end{equation}
\begin{equation}
\label{le7eq2}
\|\ba^{(t)}\|^2 \leq c_6 \cdot n < \infty,
\end{equation}
\begin{equation}
\label{le7eq3}
\min_{i=1,\dots,n, k=1,\dots,K}
\left|
\sum_{j=1}^d b_{k,j}^{(0)} \cdot x_i^{(j)} + b_{k,0}^{(0)}
\right|
\geq \delta_n >0
\end{equation}
and
\begin{equation}
\label{le7eq4}
 (d+1) \cdot t_n \cdot
\lambda_n \cdot 2 \cdot \sqrt{c_5} \cdot
\max\{1, \max_{i,l} \{|x_{i}^{(l)} |^2 \} \}
\cdot
\sqrt{c_6 \cdot n}
\cdot
\exp \left(
-
\delta_n/2
\right)
\leq
\frac{\delta_n}{2}.
\end{equation}
Then we have for
every $k \in \{1,\dots,K\}$,
any $j \in \{0,\dots,d\}$ and
any $t \in \{1,\dots,t_n\}$:
\begin{equation}
\label{le7eq5}
|b_{k,j}^{(t)}-b_{k,j}^{(t-1)}|
\leq
\lambda_n \cdot 2 \cdot \sqrt{c_5} \cdot
\max\{1, \max_{i,l} \{|x_{i}^{(l)} | \} \}
\cdot
\sqrt{c_6 \cdot n}
\cdot
\exp \left(
-
\delta_n/2
\right).
\end{equation}
\end{lemma}

\noindent
{\bf Proof.}
We show (\ref{le7eq5}) by induction on $t$. For $t=1$ the assertion follows
from
Lemma \ref{le6} and (\ref{le7eq1})-(\ref{le7eq3}). Now, we assume that  (\ref{le7eq5}) holds for all
$t \in \{1,\dots,s\}$, where $s \in \{1, \dots, t_n-1\}$. Then
\[
|b_{k,j}^{(s)}-b_{k,j}^{(0)}|
\leq
t_n \cdot
\lambda_n \cdot 2 \cdot \sqrt{c_5} \cdot
\max\{1, \max_{i,l} \{|x_{i}^{(l)} | \} \}
\cdot
\sqrt{c_6 \cdot n}
\cdot
\exp \left(
-
\delta_n/2
\right),
\]
from which, together with assumption (\ref{le7eq3}), we can conlcude that
\begin{eqnarray}
&&
\min_{i=1,\dots,n, k=1,\dots,K}
\left|
\sum_{j=1}^d b_{k,j}^{(s)} \cdot x_i^{(j)} + b_{k,0}^{(s)}
\right|
\nonumber
\\
&&
\geq
\min_{i=1,\dots,n, k=1,\dots,K}
\left|
\sum_{j=1}^d b_{k,j}^{(0)} \cdot x_i^{(j)} + b_{k,0}^{(0)}
\right|
\nonumber
\\
&&
\quad
-
\max_{i=1,\dots,n, k=1,\dots,K}
\left(
\sum_{j=1}^d |b_{k,j}^{(s)} -b_{k,j}^{(0)}| \cdot |x_i^{(j)}| + |b_{k,0}^{(s)}-b_{k,0}^{(0)}|
\right)
\nonumber
\\
&&
\geq \delta_n
- \max_{i=1,\dots,n, k=1,\dots,K}
\left(
\sum_{j=0}^d |b_{k,j}^{(s)} -b_{k,j}^{(0)}| \cdot \max\{1, \max_{i,l} \{|x_{i}^{(l)} | \} \}
\right)
\nonumber
\\
&&
\geq
\delta_n - (d+1) \cdot t_n \cdot
\lambda_n \cdot 2 \cdot \sqrt{c_5} \cdot
\max\{1, \max_{i,l} \{|x_{i}^{(l)} |^2 \} \}
\cdot
\sqrt{c_6 \cdot n}
\cdot
\exp \left(
-
\delta_n/2
\right)
\nonumber
\\
&&
\geq \frac{\delta_n}{2},
\label{ple7eq1}
\end{eqnarray}
where the last inequality is implied by inequality
 (\ref{le7eq4}).
So, for the induction step, application of Lemma \ref{le6}
together with (\ref{le7eq1}) and (\ref{ple7eq1})
yields 
\begin{eqnarray*}
	|b_{k,j}^{(s+1)}-b_{k,j}^{(s)}|
	&\leq&
	\lambda_n \cdot 2 \cdot \sqrt{F(\ba^{(s)},\bb^{(s)})} \cdot
	\max\{1, \max_{i,l} \{|x_{i}^{(l)} | \} \}
	\cdot
	|a^{(s)}_k|
	\\
	&&
	\hspace*{2cm}
	\cdot
	\exp \left(
	-
	\min_{i=1,\dots,n}
	\left\{
	\left|
	\sum_{j=1}^d b_{k,j}^{(s)} \cdot x_i^{(j)} + b_{k,0}^{(s)}
	\right|
	\right\}
	\right)
	\\
	&
	\leq &
	\lambda_n \cdot 2 \cdot \sqrt{c_5} \cdot
	\max\{1, \max_{i,l} \{|x_{i}^{(l)} | \} \}
	\cdot
	\sqrt{c_6 \cdot n}
	\cdot
	\exp \left(
	-
	\delta_n/2
	\right),
\end{eqnarray*}
from which we conclude the assertion.
\hfill $\Box$

      \begin{lemma}
  \label{le8}
Define $F$ by (\ref{se5sub2eq2}), set
\[
\lambda_n=\frac{1}{3 \cdot K}
\]
 and define $(\ba^{(t)},\bb^{(t)})$ by
(\ref{se5sub2eq3}). Assume that $(\ba^{(0)},\bb^{(0)})$ is chosen
such that
\begin{equation}
\label{le8eq1}
F(\ba^{(0)},\bb^{(0)}) \leq c_5 < \infty
\end{equation}
and
\begin{equation}
\label{le8eq2}
\min_{i=1,\dots,n, k=1,\dots,K}
\left|
\sum_{j=1}^d b_{k,j}^{(0)} \cdot x_i^{(j)} + b_{k,0}^{(0)}
\right|
\geq \delta_n \geq 1
\end{equation}
hold. Let $t_n \in \N$ and assume $2 \cdot c_1 \leq (K-2) \cdot n$,
\begin{eqnarray}
  \label{le8eq3}
  &&
  4 \cdot \max\{1,\frac{c_5}{c_1}\} \cdot \max\{1,\frac{1}{c_1^2}\} \cdot 
  \lambda_n \cdot
	(d+1)^2\cdot 
	 n^2
	\cdot
	\max\{1,\max_{i,j} |x_i^{(j)}|^4\}
	\nonumber
        \\
        &&
        \hspace*{3cm}
        \cdot
	\left(
	1+c_5 + \frac{2}{n} \sum_{i=1}^n y_i^2 
	\right)^4
	\cdot
        t_n^2
        \cdot
	\exp \left(
	-
	\delta_n/2
	\right)
	\leq 1
\end{eqnarray}
and
\begin{equation}
  \label{le8eq4}
3 \cdot t_n \cdot
\exp(-\delta_n/4) \leq 1.
\end{equation}
Then for any $t \in \{0,1,\dots,t_n-1\}$ we have 
\begin{eqnarray*}
&&
F(\ba^{(t+1)},\bb^{(t+1)})
-
\min_{\ba}
F(\ba,\bb^{(0)})
\\
&&
\leq
\left(1-
\frac{2 \cdot c_1}{3 \cdot K \cdot n}
\right)^{t+1} 
\cdot \left(
F(\ba^{(0)},\bb^{(0)})
-
\min_{\ba}
F(\ba,\bb^{(0)})
\right)
+
( 2 \sqrt{c_5} +1) 
\cdot 	\exp \left(
-
\delta_n/4
\right)
\\
&&
\quad
+ \frac{3 \cdot K \cdot n}{2 \cdot c_1} 
\cdot 
3 \cdot
\exp \left(
-
\delta_n/4
\right).
\end{eqnarray*}
\end{lemma}

      
  \noindent
      {\bf Proof.} 
We have
\begin{eqnarray*}
&&
F(\ba^{(t+1)},\bb^{(t+1)})
-
\min_{\ba}
F(\ba,\bb^{(0)})
\\
&&
=
\left(
F(\ba^{(t+1)},\bb^{(t+1)})
-
F(\ba^{(t+1)},\bb^{(t)})
\right)
+
\left(
F(\ba^{(t+1)},\bb^{(t)})
- 
\min_{\ba}
F(\ba,\bb^{(t)})
\right) \\
&&
\quad
+ 
\left(
\min_{\ba}
F(\ba,\bb^{(t)})
-
\min_{\ba}
F(\ba,\bb^{(0)})
\right)
.
\end{eqnarray*}
We will continue proving the assertion in three steps. 

\noindent {\it First step.}
We take a look at the second term
on the right-hand side of the above equality.
Lemma \ref{le2} and $|\sigma(x)| \leq 1$ give us
\begin{eqnarray*}
\left\|
(\nabla_{\ba} F)(\ba_1,\bb^{(t)})
-
(\nabla_{\ba} F)(\ba_2,\bb^{(t)})
\right\|
&&
\leq
\left(
2 \cdot (K+1) + \frac{2 \cdot c_1}{n}
\right)
\cdot \|\ba_1-\ba_2\|
\\
&&
\leq
3 \cdot K  \cdot \|\ba_1-\ba_2\|.
\end{eqnarray*}
Together with Lemma \ref{le3} this allows us to conclude from Lemma \ref{le1} that
\begin{eqnarray}
\label{le8eq**}
&&
F(\ba^{(t+1)},\bb^{(t)})
-
\min_{\ba}
F(\ba,\bb^{(t)}) \nonumber
\\
&&
\leq
\left(
1-\frac{4 \cdot c_1}{6 \cdot K \cdot n} 
\right)
\cdot
\left(
F(\ba^{(t)},\bb^{(t)})
-
\min_{\ba}
F(\ba,\bb^{(t)})
\right)
.
\end{eqnarray}
For simplicity, we introduce the following notation
\begin{eqnarray*}
	\gamma_t &=&
	\left(
	F(\ba^{(t)},\bb^{(t)})
	-
	\min_{\ba}
	F(\ba,\bb^{(0)})
	\right),
	\\
	\alpha &=& 
	\frac{2 \cdot c_1}{3 \cdot K \cdot n} 
	.
\end{eqnarray*}
As a consequence,
\begin{eqnarray}
  \label{le8eq***}
	&&
	\gamma_{t+1}
	\nonumber
	\\
	&&
	\leq 
	F(\ba^{(t+1)},\bb^{(t+1)}) - F(\ba^{(t+1)},\bb^{(t)})
	+
	(1-\alpha) \cdot (F(\ba^{(t)},\bb^{(t)}) - \min_{\ba}F(\ba,\bb^{(t)}))
	\nonumber
	\\
	&&
	\quad
	+
	\min_{\ba}F(\ba,\bb^{(t)}) - \min_{\ba}F(\ba,\bb^{(0)})
	\nonumber
	\\
	&&
	=
	F(\ba^{(t+1)},\bb^{(t+1)})-F(\ba^{(t+1)},\bb^{(t)}) 
	+ 
	(1- \alpha)\cdot (F(\ba^{(t)},\bb^{(t)})-\min_{\ba}F(\ba,\bb^{(0)}))
	\nonumber
	\\
	&&
	\quad
	+
	\alpha \cdot (\min_{\ba} F(\ba,\bb^{(t)}) - \min_{\ba} F(\ba,\bb^{(0)}))
	\nonumber
	\\
	&&
	=
	(1-\alpha) \cdot \gamma_t 
	+
	\alpha \cdot (\min_{\ba} F(\ba,\bb^{(t)}) - \min_{\ba} F(\ba,\bb^{(0)}))
	\nonumber
	\\
	&&
	\quad
	+ 
	F(\ba^{(t+1)},\bb^{(t+1)})-F(\ba^{(t+1)},\bb^{(t)}) 
	.
\end{eqnarray}

\noindent {\it Second step.} We will derive upper bounds $\beta_1,\beta_2 > 0$ such that
\begin{eqnarray*}
	\beta_1 &\geq&
	\min_{\ba}
	F(\ba,\bb^{(t)})
	-
	\min_{\ba}
	F(\ba,\bb^{(0)})
	,
	\\
	\beta_2 &\geq&
	F(\ba^{(t+1)},\bb^{(t+1)})
	-
	F(\ba^{(t+1)},\bb^{(t)}).
\end{eqnarray*}

\noindent We will start with finding $\beta_1$.
%
%
In the process we will also derive an upper bound $\beta_2$. 

\noindent
Choose $\bar{\ba}$ such that
\[
F(\bar{\ba},\bb^{(0)})
=
\min_{\ba}
F(\ba,\bb^{(0)}).
\]
Then
\[
\frac{c_1}{n} \cdot \sum_{k=0}^n \bar{\ba}_k^2
\leq
F(\bar{\ba},\bb^{(0)})
\leq
F(\ba^{(0)},\bb^{(0)})
\leq c_5,
\]
hence
\[
\sum_{k=0}^K \bar{a}_k^2
\leq \frac{c_5 \cdot n}{c_1}.
\]
We have
\begin{eqnarray*}
&&
\min_{\ba}
F(\ba,\bb^{(t)})
-
\min_{\ba}
F(\ba,\bb^{(0)})
=
\min_{\ba}
F(\ba,\bb^{(t)})
-
F(\bar{\ba},\bb^{(0)})
\\
&&
\leq
F(\bar{\ba},\bb^{(t)})
-
F(\bar{\ba},\bb^{(0)})
\\
&&
=
\frac{1}{n} \sum_{i=1}^n
( f_{net,(\bar{\ba},\bb^{(t)})}(x_i)
+f_{net,(\bar{\ba},\bb^{(0)})}(x_i)
- 2 y_i)
\cdot
(f_{net,(\bar{\ba},\bb^{(t)})}(x_i)
-f_{net,(\bar{\ba},\bb^{(0)})}(x_i)
)
\\
&&
=
\frac{1}{n} \sum_{i=1}^n
 ( 2f_{net,(\bar{\ba},\bb^{(0)})}(x_i)
- 2 y_i)
\cdot
(f_{net,(\bar{\ba},\bb^{(t)})}(x_i)
-f_{net,(\bar{\ba},\bb^{(0)})}(x_i)
)
\\
&&
\quad
+
\frac{1}{n} \sum_{i=1}^n
(f_{net,(\bar{\ba},\bb^{(t)})}(x_i)
-f_{net,(\bar{\ba},\bb^{(0)})}(x_i)
)
^2
\\
&&
\leq
2
\cdot
\sqrt{F(\bar{\ba},\bb^{(0)})}
\cdot
\sqrt{
\frac{1}{n} \sum_{i=1}^n
(f_{net,(\bar{\ba},\bb^{(t)})}(x_i)
-f_{net,(\bar{\ba},\bb^{(0)})}(x_i)
)
^2
}
\\
&&
\quad
+
\frac{1}{n} \sum_{i=1}^n
(f_{net,(\bar{\ba},\bb^{(t)})}(x_i)
-f_{net,(\bar{\ba},\bb^{(0)})}(x_i)
)
^2.
\end{eqnarray*}
Applying the Cauchy-Schwarz inequality a second time 
and since $\sigma$ is Lipschitz continuous
we get
\begin{eqnarray*}
&&
\frac{1}{n} \sum_{i=1}^n
 (f_{net,(\bar{\ba},\bb^{(t)})}(x_i)
-f_{net,(\bar{\ba},\bb^{(0)})}(x_i)
)
^2
\\
&&
=
\frac{1}{n} \sum_{i=1}^n
\left(
\sum_{k=1}^K
\bar{a}_k
\cdot
\left(
\sigma
\left(
\sum_{j=1}^d b_{k,j}^{(t)} \cdot x_i^{(j)} + b_{k,0}^{(t)}
\right)
-
\sigma
\left(
\sum_{j=1}^d b_{k,j}^{(0)} \cdot x_i^{(j)} + b_{k,0}^{(0)}
\right)
\right)
\right)^2
\\
&&
\leq
\sum_{k=1}^K
\bar{a}_k^2
\cdot
\max\{1,\max_{i,j} |x_i^{(j)}|^2\}
\cdot
(d+1)
\cdot
\sum_{k=1}^K
\sum_{j=0}^d| b_{k,j}^{(t)} - b_{k,j}^{(0)}|^2.
\end{eqnarray*}
By Lemma \ref{le7} 
where, as we will show below, $c_5$ and $c_6$ are replaced by
\[
1+c_5 + \frac{2}{n} \sum_{i=1}^n y_i^2
\quad \mbox{and} \quad
\left( 1+c_5 + \frac{2}{n} \sum_{i=1}^n y_i^2 \right)
\cdot
\frac{1}{c_1}, \mbox{ respectively,}
\]
we know that for any $k \in\{1,\dots,K\}$ and any $j \in \{0,\dots,d\}$
\begin{eqnarray*}
	&&
	| b_{k,j}^{(t)} - b_{k,j}^{(0)}|
	\\
	&&
	\leq
	| b_{k,j}^{(t)} - b_{k,j}^{(t-1)}| + | b_{k,j}^{(t-1)} - b_{k,j}^{(t-2)}| + \cdots + | b_{k,j}^{(1)} - b_{k,j}^{(0)}|	
	\\
	&&
	\leq
	t \cdot 
	\lambda_n \cdot 2 \cdot  \left( 1+ c_5 + \frac{2}{n} \sum_{i=1}^n y_i^2 \right) \cdot
        \max\{1,\frac{1}{c_1}\} \cdot
	\max\{1, \max_{i,l} \{|x_{i}^{(l)} | \} \}
	\cdot
	\sqrt{n}
	\cdot
	\exp \left(
	-
	\delta_n/2
	\right).
\end{eqnarray*}
From this we conclude that
\begin{eqnarray*}
	&&
	\frac{1}{n} \sum_{i=1}^n
	(f_{net,(\bar{\ba},\bb^{(t)})}(x_i)
	-f_{net,(\bar{\ba},\bb^{(0)})}(x_i)
	)
	^2
	\\
	&&
	\leq
	\sum_{k=1}^K
	\bar{a}_k^2
	\cdot
	\max\{1,\max_{i,j} |x_i^{(j)}|^2\}
	\cdot
	(d+1)
	\cdot
	K \cdot
	(d+1)  
	\\
	&&
	\quad
	\cdot  \left(
		t \cdot 
	        \lambda_n \cdot \left( 1+ c_5 + \frac{2}{n} \sum_{i=1}^n y_i^2 \right)
                \cdot
                \max\{1,\frac{1}{c_1}\}
                \cdot
	\max\{1, \max_{i,l} \{|x_{i}^{(l)} | \} \}
	\cdot
	\sqrt{n}
	\cdot
	\exp \left(
	-
	\delta_n/2
	\right)
	\right)^2
	\\
	&&
	\leq
	t^2 \cdot
        \frac{c_5}{c_1} \cdot
        \max\{1,\frac{1}{c_1^2} \} \cdot
	\left( 1+ c_5 + \frac{2}{n} \sum_{i=1}^n y_i^2 \right)^2 \cdot n^2
		\cdot
	\max\{1,\max_{i,j} |x_i^{(j)}|^4\}
	\cdot
	(d+1)^2
        \\
        &&
        \hspace*{5cm}
	\cdot
	K \cdot
	\lambda_n^2\cdot  
	\exp \left(
	-
	\delta_n
	\right)
        \\
        &&
        \leq \exp \left(
	-
	\delta_n/2
	\right),
\end{eqnarray*}
where the last inequality follows from (\ref{le8eq3}).
Hence,
\begin{eqnarray*}
  &&
  \min_{\ba}
F(\ba,\bb^{(t)})
-
\min_{\ba}
F(\ba,\bb^{(0)})
\\
	&&
	\leq
	2
	\cdot
	\sqrt{F(\bar{\ba},\bb^{(0)})}
	\cdot
	\sqrt{
		\frac{1}{n} \sum_{i=1}^n
		(f_{net,(\bar{\ba},\bb^{(t)})}(x_i)
		-f_{net,(\bar{\ba},\bb^{(0)})}(x_i)
		)
		^2
	}
	\\
	&&
	\quad
	+
	\frac{1}{n} \sum_{i=1}^n
	(f_{net,(\bar{\ba},\bb^{(t)})}(x_i)
	-f_{net,(\bar{\ba},\bb^{(0)})}(x_i)
	)
	^2
        \\
	&&
	\leq
	(2 \cdot \sqrt{F(\bar{\ba},\bb^{(0)})}+1)
	\cdot
	\sqrt{
		\frac{1}{n} \sum_{i=1}^n
		(f_{net,(\bar{\ba},\bb^{(t)})}(x_i)
		-f_{net,(\bar{\ba},\bb^{(0)})}(x_i)
		)
		^2
	}
	\\
	&&
	\leq
        ( 2 \cdot \sqrt{c_5} +1) \cdot
\exp \left(
	-
	\delta_n/4
	\right)
        =
	\beta_1.
\end{eqnarray*}
\noindent It remains to be shown that Lemma $\ref{le7}$ was, in fact, applicable, i.e. we will show that the conditions of Lemma $\ref{le7}$ are met. 
For that we show the following claim for all
$s \in \{0,1,\dots,t_n-1\}$
by induction
\begin{eqnarray}
&&
\label{le8eq*}
\max\left\{
F(\ba^{(s+1)},\bb^{(s)}),
F(\ba^{(s+1)},\bb^{(s+1)})
\right\}
-
\min_{\ba}
F(\ba,\bb^{(0)}) \nonumber \\
&&
\leq
c_5 + \frac{1}{n} \sum_{i=1}^n y_i^2 + 3 \cdot (s+1) \cdot \exp(-\delta_n/4).
\end{eqnarray}
While doing so, we will be deriving an upper bound $\beta_2$ in the process.
For $s=0$ the inequality trivially holds by
(\ref{le8eq**}), (\ref{le8eq***}),
(\ref{le8eq1}) and by the bound
\[
F(\ba^{(1)},\bb^{(1)})
-
F(\ba^{(1)},\bb^{(0)})
\leq
3 \cdot
\exp \left(
-
\delta_n/4
\right)
\]
which will be proven below (cf., (\ref{le8eq****})).

So, for the induction hypothesis, assume that (\ref{le8eq*}) holds for
$s=t-1$
for arbitrary $t \in \{1, \dots, t_n-1\}$.
Trivially we have
\[
\min_{\ba}
F(\ba,\bb^{(t)})
-
\min_{\ba}
F(\ba,\bb^{(0)})
\leq
F(\mathbf{0},\bb^{(t)})
=
\frac{1}{n} \sum_{i=1}^n y_i^2,
\]
hence by (\ref{le8eq**}) and by the induction assumption we get
\begin{eqnarray}
&&
F(\ba^{(t+1)},\bb^{(t)})
-
\min_{\ba}
F(\ba,\bb^{(0)})
\nonumber \\
&&
\leq
\left(
1-\frac{2 \cdot c_1}{3 \cdot K \cdot n} 
\right)
\cdot
\left(
F(\ba^{(t)},\bb^{(t)})
-
\min_{\ba}
F(\ba,\bb^{(0)})
\right)
\nonumber \\
&&
\quad
+
\frac{2\cdot c_1}{3 \cdot K \cdot n}
\cdot
\left(
\min_{\ba}
F(\ba,\bb^{(t)})
-
\min_{\ba}
F(\ba,\bb^{(0)})
\right)
\nonumber \\
&&
\leq
\left(
1-\frac{2\cdot c_1}{3 \cdot K \cdot n} 
\right)
\cdot
\left( c_5 + \frac{1}{n} \sum_{i=1}^n y_i^2
+ 3 \cdot t \cdot \exp(-\delta_n/4)
\right)
\nonumber \\
&&
\quad
+
\frac{2 \cdot c_1}{3 \cdot K \cdot n}
\cdot
\frac{1}{n} \sum_{i=1}^n y_i^2
\nonumber \\
&&
\leq
c_5 + \frac{1}{n} \sum_{i=1}^n y_i^2+ 3 \cdot t \cdot \exp(-\delta_n/4).
\label{le8eq30}
\end{eqnarray}
Next, by (\ref{le8eq***}) and by the induction hypothesis
we get
\begin{eqnarray}
&&
F(\ba^{(t+1)},\bb^{(t+1)})
-
\min_{\ba}
F(\ba,\bb^{(0)})
\nonumber \\
&&
\leq
c_5 + \frac{1}{n} \sum_{i=1}^n y_i^2+ 3 \cdot t \cdot \exp(-\delta_n/4)
+
F(\ba^{(t+1)},\bb^{(t+1)})
-
F(\ba^{(t+1)},\bb^{(t)}).
\label{le8eq5*}
\end{eqnarray}
Further, we have
\begin{eqnarray*}
	&&
	F(\ba^{(t+1)},\bb^{(t+1)})
	-
	F(\ba^{(t+1)},\bb^{(t)})
	\\
	&&
	=
	\frac{1}{n} \sum_{i=1}^n
	( f_{net,(\ba^{(t+1)},\bb^{(t+1)})}(x_i)
	+f_{net,\ba^{(t+1)},\bb^{(t)})}(x_i)
	- 2 y_i)
	\\
	&&
	\hspace*{1cm}
	\cdot
	(f_{net,(\ba^{(t+1)},\bb^{(t+1)})}(x_i)
	-f_{net,(\ba^{(t+1)},\bb^{(t)})}(x_i)
	)
	\\
	&&
	=
	\frac{1}{n} \sum_{i=1}^n
	( 2f_{net,(\ba^{(t+1)},\bb^{(t)})}(x_i)
	- 2 y_i)
	\cdot
	(f_{net,(\ba^{(t+1)},\bb^{(t+1)})}(x_i)
	-f_{net,(\ba^{(t+1)},\bb^{(t)})}(x_i)
	)
	\\
	&&
	\quad
	+
	\frac{1}{n} \sum_{i=1}^n
	(f_{net,(\ba^{(t+1)},\bb^{(t+1)})}(x_i)
	-f_{net,(\ba^{(t+1)},\bb^{(t)})}(x_i)
	)
	^2
	\\
	&&
	\leq
	2
	\cdot
	\sqrt{F(\ba^{(t+1)},\bb^{(t)})}
	\cdot
	\sqrt{
		\frac{1}{n} \sum_{i=1}^n
		(f_{net,(\ba^{(t+1)},\bb^{(t+1)})}(x_i)
		-f_{net,(\ba^{(t+1)},\bb^{(t)})}(x_i)
		)
		^2
	}
	\\
	&&
	\quad
	+
	\frac{1}{n} \sum_{i=1}^n
	(f_{net,(\ba^{(t+1)},\bb^{(t+1)})}(x_i)
	-f_{net,(\ba^{(t+1)},\bb^{(t)})}(x_i)
	)
	^2.
\end{eqnarray*}
Since $\sigma$ is Lipschitz continuous, applying the Cauchy-Schwarz inequality once more
yields
\begin{eqnarray*}
	&&
	\frac{1}{n} \sum_{i=1}^n
	(f_{net,(\ba^{(t+1)},\bb^{(t+1)})}(x_i)
	-f_{net,(\ba^{(t+1)},\bb^{(t)})}(x_i)
	)
	^2
	\\
	&&
	=
	\frac{1}{n} \sum_{i=1}^n
	\left(
	\sum_{k=1}^K
	(\ba^{(t+1)})_k
	\cdot
	\left(
	\sigma
	\left(
	\sum_{j=1}^d b_{k,j}^{(t+1)} \cdot x_i^{(j)} + b_{k,0}^{(t+1)}
	\right)
	-
	\sigma
	\left(
	\sum_{j=1}^d b_{k,j}^{(t)} \cdot x_i^{(j)} + b_{k,0}^{(t)}
	\right)
	\right)
	\right)^2
	\\
	&&
	\leq
	\sum_{k=1}^K
	(\ba^{(t+1)})_k^2
	\cdot
	\max\{1,\max_{i,j} |x_i^{(j)}|^2\}
	\cdot
	(d+1)
	\cdot
	\sum_{k=1}^K
	\sum_{j=0}^d| b_{k,j}^{(t+1)} - b_{k,j}^{(t)}|^2 
	\\
	&&
	\leq
	\frac{n}{c_1} \cdot F(\ba^{(t+1)},\bb^{(t)})
	\cdot
	\max\{1,\max_{i,j} |x_i^{(j)}|^2\}
	\cdot
	(d+1)
	\cdot
	\sum_{k=1}^K
	\sum_{j=0}^d| b_{k,j}^{(t+1)} - b_{k,j}^{(t)}|^2 
	.
\end{eqnarray*}
By Lemma \ref{le7} (where (\ref{le7eq1}) and (\ref{le7eq2}) are true because
of the fact that
the induction hypothesis implies that we have
\[
F(\ba^{(t)},\bb^{(t)})
\leq
c_5 + \frac{1}{n} \sum_{i=1}^n y_i^2 + 1 + F(\mathbf{0},\bb^{(0)})
\leq
1+c_5 + \frac{2}{n} \sum_{i=1}^n y_i^2
,
\]
from which (together with the defnition of $F$) we can conclude
that (\ref{le7eq1}) and (\ref{le7eq2}) hold if we replace
there $c_5$ and $c_6$ by
\[
1+c_5 + \frac{2}{n} \sum_{i=1}^n y_i^2
\quad \mbox{and} \quad
\left( 1+c_5 + \frac{2}{n} \sum_{i=1}^n y_i^2 \right)
\cdot
\frac{1}{c_1}, \mbox{ respectively,}
\]
and where (\ref{le7eq4}) holds because of (\ref{le8eq3}))
and because of (\ref{le8eq2})
we know that for any $k \in\{1,\dots,K\}$ and any $j \in \{0,\dots,d\}$
we have
\begin{eqnarray*}
	&&
	| b_{k,j}^{(t+1)} - b_{k,j}^{(t)}|
	\\
	&&
	\leq 
	\lambda_n \cdot 2 \cdot
	\left( 1+c_5 + \frac{2}{n} \sum_{i=1}^n y_i^2 \right)
	\cdot
	\max\{1, \max_{i,l} \{|x_{i}^{(l)} | \} \}
	\cdot
	\sqrt{n}
	\cdot
	\exp \left(
	-
	\delta_n/2
	\right) / \sqrt{c_1}.
\end{eqnarray*}
Together with
\begin{eqnarray*}
	F(\ba^{(t+1)},\bb^{(t)})
	&\leq&
	\min_{\ba} F(\ba,\bb^{(0)})
	+ c_5 + \frac{1}{n} \sum_{i=1}^n y_i^2 + 3 \cdot t \cdot \exp(-\delta_n/4)
	\\
	&\leq&
	1+c_5 + \frac{2}{n} \sum_{i=1}^n y_i^2,
\end{eqnarray*}
where the first inequality follows trivially from $(\ref{le8eq30})$,
this implies
\begin{eqnarray*}
	&&
	\frac{1}{n} \sum_{i=1}^n
	(f_{net,(\ba^{(t+1)},\bb^{(t+1)})}(x_i)
	-f_{net,(\ba^{(t+1)},\bb^{(t)})}(x_i)
	)
	^2
	\\
	&&
	\leq
	4 \cdot
	\frac{n^{2}}{c_1^2} \cdot \lambda_n^2
	\cdot
	\left( 1+c_5 + \frac{2}{n} \sum_{i=1}^n y_i^2 \right)^3
	\cdot
	\max\{1,\max_{i,j} |x_i^{(j)}|^4\}
	\cdot
	(d+1)^2 \cdot K
	\cdot
	\exp \left(
	-
	\delta_n
	\right)
	\\
	&&
	\leq
	4 \cdot
	\frac{(d+1)^2 \cdot n^2}{c_1^2}
	\cdot
	\max\{1,\max_{i,j} |x_i^{(j)}|^4\}
	\cdot
	\left(
	1+c_5 + \frac{2}{n} \sum_{i=1}^n y_i^2 
	\right)^4
	\cdot
	\exp \left(
	-
	\delta_n/2
	\right)
	\\
	&&
	\quad \cdot
	\min \left\{1,
	(F(\ba^{(t+1)},\bb^{(t)}))^{-1}
	\right\} \cdot
	\exp \left(
	-
	\delta_n/2
	\right)        \\
	&&
	\leq
	\min \left\{1,
	(F(\ba^{(t+1)},\bb^{(t)}))^{-1}
	\right\} \cdot
	\exp \left(
	-
	\delta_n/2
	\right). 
\end{eqnarray*}
(Here, the last inequality follows from (\ref{le8eq3}).)
Summarizing the above results we get
\begin{equation}
\label{le8eq****}
F(\ba^{(t+1)},\bb^{(t+1)})
-
F(\ba^{(t+1)},\bb^{(t)})
\leq
3 \cdot
\exp \left(
-
\delta_n/4
\right)
= 
\beta_2.
\end{equation}
By combining this inequality with the results
above we get (\ref{le8eq*}) for $s=t$. This concludes the proof of
(\ref{le8eq*}).
 Thus, all the conditions of Lemma \ref{le7} are met, since we can conclude from 
\[
\min_{\ba} F(\ba,\bb^{(0)})
\leq
F(\mathbf{0},\bb^{(0)})
=
\frac{1}{n} \sum_{i=1}^n y_i^2
\]
and from inequalities (\ref{le8eq*}) and  (\ref{le8eq4})
that also
(\ref{le7eq1})
and (because of the defintion of $F$) (\ref{le7eq2})
hold where $c_5$ and $c_6$ are replaced by
\[
1+c_5 + \frac{2}{n} \sum_{i=1}^n y_i^2
\quad \mbox{and} \quad
\left( 1+c_5 + \frac{2}{n} \sum_{i=1}^n y_i^2 \right)
\cdot
\frac{1}{c_1}, \mbox{ respectively.}
\]
As above we also see that (\ref{le7eq4}) holds.

{\it Third Step.} The results we derived in the first step imply that
\[
\gamma_{t+1} \leq (1-\alpha) \cdot \gamma_t + \alpha \cdot \beta_1 + \beta_2.
\]
Applying this relation recursively 
using standard techniques from the literature
we get
\begin{eqnarray*}
	\gamma_{t+1}
	&\leq&
	(1-\alpha)\cdot((1 - \alpha) \cdot \gamma_{t-1} + \alpha \cdot \beta_1 + \beta_2) + \alpha \cdot \beta_1 + \beta_2
	\\
	&=&
	(1-\alpha)^2 \cdot \gamma_{t-1} + (1-\alpha) \cdot \alpha \cdot \beta_1
	+ \alpha \cdot    \beta_1 + (1-\alpha)\cdot \beta_2 + \beta_2
	\\
	& \leq & \dots \\
	&\leq&
	(1-\alpha)^{t+1} \cdot \gamma_{0} + \sum_{k = 0}^{t} (1-\alpha)^k \cdot \alpha \cdot \beta_1 + \sum_{k = 0}^{t} (1-\alpha)^k \cdot \beta_2
	\\
	&\leq&
	(1-\alpha)^{t+1} \cdot \gamma_{0} + \sum_{k = 0}^{\infty} (1-\alpha)^k \cdot \alpha \cdot \beta_1 + \sum_{k = 0}^{\infty} (1-\alpha)^k \cdot \beta_2
	\\
	&=&
	(1-\alpha)^{t+1} \cdot \gamma_{0} + \frac{\alpha \cdot \beta_1}{1 - (1 - \alpha)} + \frac{ \beta_2}{1 - (1 - \alpha)} 
	\\
	&=&
	(1-\alpha)^{t+1} \cdot \gamma_{0} +\beta_1 + \frac{\beta_2}{\alpha}.
\end{eqnarray*}
Plugging in the above results yields
\begin{eqnarray*}
	&&
	F(\ba^{(t+1)},\bb^{(t+1)})
	-
	\min_{\ba}
	F(\ba,\bb^{(0)})
	\\
	&&
	\leq
	(1-\alpha)^{t+1} \cdot \gamma_{0} +\beta_1 + \frac{\beta_2}{\alpha}.
	\\
	&&
	\leq
	\left(1-
	\frac{2\cdot c_1}{3 \cdot K  \cdot n} 
	\right)^{t+1} 
	\cdot \left(
	F(\ba^{(0)},\bb^{(0)})
	-
	\min_{\ba}
	F(\ba,\bb^{(0)})
	\right)
        +
             ( 2 \cdot \sqrt{c_5} +1) \cdot
\exp \left(
	-
	\delta_n/4
	\right)
   	\\
	&&
	\quad
	+ \frac{3 \cdot K \cdot n}{2 \cdot c_1} 
	\cdot 3 \cdot
	\exp \left(
	-
	\delta_n/4
	\right),
\end{eqnarray*}
which concludes the proof.

      \hfill $\Box$

\subsection{Two auxiliary results from empirical process theory}
\label{se5sub2}
\begin{lemma}
\label{le9}
Let $\beta_n = c_{3} \cdot \log(n)$ for some suitably large
constant $c_{3}>0$. Assume that the distribution of $(X,Y)$ satisfies
(\ref{subgaus})
for some constant $c_{2}>0$ and that the regression function $m$ is
bounded in absolute value. Let
$\F_n$ be a set of functions $f:\Rd \rightarrow \R$ and assume that
the estimate $m_n$ satisfies
\[
m_n=T_{\beta_n}\tilde{m}_n,
\]
\[
\tilde{m}_n(\cdot)
=
\tilde{m}_n(\cdot,(X_1,Y_1),\dots,(X_n,Y_n))
\in \F_n
\]
and
\begin{eqnarray*}
  &&
\frac{1}{n}
\sum_{i=1}^n |Y_i-\tilde{m}_n(X_i)|^2
\cdot
I_{
\{ |Y_i| \leq \beta_n \mbox{ for all } i \in \{1, \dots, n\}\}
  }
\\
&&
\leq
\min_{l \in \Theta_n}
\left(
\frac{1}{n}
\sum_{i=1}^n |Y_i-g_{n,l}(X_i)|^2+pen_n(g_{n,l}) + \epsilon_{n,l}
\right)
\end{eqnarray*}
for some random functions $g_{n,l} : \Rd \rightarrow \R$, some
nonempty parameter set $\Theta_n$ and some 
deterministic penalty terms $pen_n(g_{n,l}) \geq 0$, and some
additional deterministic term $\epsilon_{n,l}$,
where the functions $g_{n,l}$ only depend on the set
\[
\mathcal{D}_{n,r} = \{X_1,\dots,X_n, \bar{\bc}_1^{(1)}, \dots, \bar{\bc}_r^{(1)}, \dots, \bar{\bc}_1^{(I_n)}, \dots, \bar{\bc}_r^{(I_n)}\}
\]
and where $\bar{\bc}_1^{(1)}, \dots, \bar{\bc}_r^{(1)}, \dots, \bar{\bc}_1^{(I_n)}$ are random variables independent of $(X_1,Y_1), \dots, (X_n, Y_n)$.

Then $m_n$ satisfies
\begin{align*}
 & \mathbf E \int |m_n(x) - m(x)|^2 \PROB_X (dx) \leq
                                                  \frac{c_{7}\cdot
                                                  (\log n)^2\cdot
    \left(
    \log\left(
\sup_{x_1^n}
\mathcal{N}_1 \left(\frac{1}{n\cdot\beta_n}, \F_n, x_1^n \right)
\right)
+1
\right)
  }{n}\\
                                                  & \hspace*{4.2cm}
                                                +  2 \cdot  
                                                \EXP
                                                \left(
\min_{l \in \Theta_n}
\frac{1}{n} \sum_{i = 1}^{n}
|g_{n,l}(X_i)-m(X_i)|^2
+ pen_n(g_{n,l})
+ {\epsilon}_{n,l}
\right)
\end{align*}
for $n>1$ and some constant $c_{7}>0$, which does not depend on $n, \beta_n$ or the parameters of the estimate.
\end{lemma}
\noindent
    {\bf Proof.}
    This lemma follows in a straightforward way from the proof of
    Theorem 1 in Bagirov, Clausen and Kohler (2009). A complete version of the proof
    is given in the Supplement. 
    \hfill $\Box$\\

\noindent
In order to bound the covering number
$\mathcal{N}_1 \left(\frac{1}{n\cdot\beta_n}, \mathcal{F}_n,
  x_1^n \right)$
we will use the following lemma.

\begin{lemma}
\label{le10}
Let $\max\{K,\beta_n,\gamma_n\} \leq n^{c_8}$ and define $\F$ by
\begin{eqnarray*}
\F&=& \Bigg\{
f:\Rd \rightarrow \R \, : \,
f(x)=\sum_{k=0}^K a_k \cdot \sigma \left(
\sum_{j=1}^d b_{k,j} \cdot x^{(j)} + b_{k,0}
\right)
\quad (x \in \Rd) 
\\
&&
\hspace*{3cm}
\mbox{ for some }
a_k,b_{k,j} \in \R
\mbox{ satisfying }
\sum_{k=0}^K a_k^2 \leq \gamma_n.
\Bigg\}
\end{eqnarray*}
Then we have for any $x_1^n \in (\Rd)^n$:
\[
\log\left(
\mathcal{N}_1 \left(\frac{1}{n\cdot\beta_n}, \mathcal{F},
 x_{1}^n\right) \right)
\leq
c_{9} \cdot \log n \cdot K.
\]

\end{lemma}

\noindent
    {\bf Proof.}
Using that
\[
\sum_{k=0}^{K} |a_k|^2
\leq \gamma_n
\]
implies
\[
\sum_{k=0}^{K} |a_k|
\leq
\sqrt{K+1} \cdot \sqrt{
\sum_{k=0}^{K} |a_k|^2
}
\leq
\sqrt{(K+1) \cdot \gamma_n},
\]
we can conclude from Lemma 16.6 in  Gy\"orfi et al. (2002)
that we have
\begin{eqnarray*}
  &&
  \Nu_1 \left(
\frac{1}{n \cdot \beta_n},\F,x_1^n
\right)
\\
&&
\leq
\left(
\frac{e ( \sqrt{(K+1) \gamma_n} + 1/(n \cdot \beta_n))}
     {1/(2 \cdot n \cdot \beta_n)}
     \right)^{K+1}
     \cdot
     \left(
     \Nu_1 \left(
     \frac{1/(2 \cdot n \cdot \beta_n)}{\sqrt{(K+1) \gamma_n} + 1/(n \cdot \beta_n)},
     \G, x_1^n
     \right)
     \right)^{K+1},
  \end{eqnarray*}
where
\begin{eqnarray*}
  \G &=&
  \Bigg\{
g:\Rd \rightarrow \R \, : \,
g(x)=\sigma \left(
\sum_{j=1}^d b_{j} \cdot x^{(j)} + b_{0}
\right)
\quad (x \in \Rd) 
\\
&&
\hspace*{6cm}
\mbox{ for some }
b_0, \dots, b_d \in \R
\Bigg\}.
  \end{eqnarray*}
By Lemma 16.3, Theorem 9.5 and
Theorem 9.4 in  Gy\"orfi et al. (2002) we get
\begin{eqnarray*}
  &&
     \Nu_1 \left(
     \frac{1/(2 \cdot n \cdot \beta_n)}{\sqrt{(K+1) \gamma_n} + 1/(n \cdot \beta_n)},
     \G, x_1^n
     \right)
     \\
     &&
     \leq
     3
     \cdot
     \Bigg(
     2e \cdot
     (2 \cdot n \cdot \sqrt{K+1} \cdot \beta_n \cdot \sqrt{\gamma_n}
     +
     2)
     \\
     &&
     \hspace*{4cm}
     \cdot
     \log
     \left(
     3e \cdot
     (2 \cdot n \cdot \sqrt{K+1} \cdot \beta_n \cdot \sqrt{\gamma_n}
     +
     2)
     \right)
     \Bigg)^{d+2},
  \end{eqnarray*}
which implies the assertion.
 \quad \hfill $\Box$

\subsection{Proof of Theorem \ref{th1}}
On the event
\[
B_n = \{ |Y_i| \leq \sqrt{n} \, : \, i=1, \dots, n\}
\]
we know by (\ref{se2eq2}) that we have $\tilde{m}_n \in \F$, where $\F$ is the function set defined in Lemma \ref{le10} and where we set $\gamma_n = \sqrt{n}$. 
Define the estimate $\bar{m}_n$ by
\[
\bar{m}_n = 
\begin{cases} 
m_n  &\mbox{if } B_n \\ 
0 & \mbox{if } B_n^c.
\end{cases}
\]
Then,
\[
\int |m_n(x)-m(x)|^2 \PROB_X(dx)
\leq
\int |\bar{m}_n(x)-m(x)|^2 \PROB_X(dx) + 2 \beta_n^2 \cdot 1_{B_n^c}.
\]
By Markov inequality we know
\[
\PROB\{B_n^c\} \leq n \cdot \PROB\{|Y|> \sqrt{n}\}
\leq
\frac{n \cdot \EXP\{ e^{c_3 \cdot Y^2} \} }{\exp(c_3 \cdot n)},
\]
therefore (\ref{subgaus}) implies that it suffices to show the assertion
under the additional assumption
\begin{equation*}
\label{pth2eq1}
\tilde{m}_n(\cdot, (X_1,Y_1), \dots, (X_n,Y_n))
\in \F.
\end{equation*}

From the definition of the estimate and from Lemma \ref{le8}
we get
\begin{eqnarray*}
  &&
\frac{1}{n} \sum_{i=1}^n |Y_i-\tilde{m}_n(X_i)|^2
\cdot
I_{
\{ |Y_i| \leq \beta_n \mbox{ for all } i \in \{1, \dots, n\}\}
}
\\
&&
\leq
\min_{\ba \in \R^{K+1},l=1,\dots,I_n}
\left(
\frac{1}{n} \sum_{i=1}^n |Y_i-f_{net,(\ba,(\bb^{(l)})^{(0)})} (X_i)|^2
+
\frac{c_1}{n} \sum_{k=0}^{K \cdot r} a_k^2 + {\epsilon}_n
\right)
\end{eqnarray*}
where
\begin{eqnarray*}
  {\epsilon}_n
  &=&
  \left(
1 - \frac{2 \cdot c_1}{3 \cdot K \cdot n}
\right)^{t_n}
\cdot
\beta_n^2
+
(2 \cdot \beta_n+1) \cdot
      \exp \left(
- \frac{\sqrt{d} \cdot A \cdot \rho_n}{4\cdot (n+1) \cdot (K-1)}
\right)
\\
&&
\hspace*{3cm}
+
\frac{3 \cdot K \cdot n}{2 \cdot c_1}
\cdot 
3 \cdot
\exp \left(
- \frac{\sqrt{d} \cdot A \cdot \rho_n}{4 (n+1) \cdot (K-1)}
\right)
\\
&\leq&
\exp \left( - \frac{2 \cdot c_1}{3} \cdot (\log n)^2 \right)
\cdot
\beta_n^2
+
( 2 \cdot \beta_n +1)
\cdot \exp \left(
- \frac{\sqrt{d} \cdot A \cdot n}{8}
\right) 
\\
&&
\hspace*{3cm}
+
\frac{3 \cdot K \cdot n}{2 \cdot c_1}
\cdot 
3 \cdot
\exp \left(
- \frac{\sqrt{d} \cdot A}{8} \cdot n
\right) .
\end{eqnarray*}
Application of Lemma \ref{le9} and of Lemma \ref{le10} yields
\begin{eqnarray*}
&&
\EXP \int |m_n(x)-m(x)|^2 \PROB_X (dx)
\\
&&
\leq
c_{10} \cdot \frac{(\log n)^3 \cdot K \cdot r}{n}
\\
&&
\quad
+
2
\cdot 
\EXP \left(
 \min_{\ba \in \R^{K+1},l=1,\dots,I_n}
 \frac{1}{n} \sum_{i = 1}^{n}
 | f_{net,(\ba,(\bb^{(l)})^{(0)})}  (X_i)-m(X_i)|^2
+
\frac{c_1}{n} \sum_{k=0}^{K \cdot r} a_k^2 
\right)
\\
&&
\quad
+
2 \cdot \exp \left( - \frac{2 \cdot c_1}{3} \cdot (\log n)^2 \right)
\cdot
\beta_n^2
+
2 \cdot
(2 \cdot \beta_n+1) \cdot \exp \left(
- \frac{\sqrt{d} \cdot A  \cdot n}{8}
\right)
\\
&&
\quad
+ 
2
\cdot
\frac{3 \cdot K \cdot n}{2 \cdot c_1}
\cdot 
3 \cdot
\exp \left(
- \frac{\sqrt{d} \cdot A}{8} \cdot n
\right).
\end{eqnarray*}
We have for all $x \in [-A,A]^d$
\begin{eqnarray*}
	&&
| f_{net,(\ba,(\bb^{(l)})^{(0)})}  (x)-m(x)|
\\
&&
\leq
|m(x) - \sum_{s=1}^{r} g_s((\bar{\bc}^{(l)})^T_s x)| 
+
| \sum_{s=1}^{r} g_s((\bar{\bc}^{(l)})^T_s x) - f_{net,(\ba,(\bb^{(l)})^{(0)})}  (x) |
\end{eqnarray*}
The $(p,C)$-smoothness of the $g_s$ implies for all $x \in [-A,A]^d$
\begin{eqnarray*}
|m(x) -\sum_{s=1}^r g_s((\bar{\bc}^{(l)})^T_s x)|
&=&
|\sum_{s=1}^r g_s (\bc_s^T x) -\sum_{s=1}^r g_s((\bar{\bc}^{(l)})^T_s x)|
\\
&
\leq &
\sum_{s=1}^r C \cdot |\bc_s^T x - (\bar{\bc}^{(l)})^T_s x|^p
\\
&
\leq &
r \cdot C \cdot \sup_{x \in [-A,A]^d}\|x\|^p \cdot \max_{s=1,\dots,r}\|\bc_s - \bar{\bc}^{(l)}_s \|_{\infty}^p.
\end{eqnarray*}
By Lemma $\ref{le5}$ we get for all $x \in [-A,A]^d$
\begin{eqnarray*}
	&&
	| \sum_{s=1}^{r} g_s((\bar{\bc}^{(l)})^T_s x) - f_{net,(\ba,(\bb^{(l)})^{(0)})}  (x) |
	\\
	&&
	=
	| \sum_{s=1}^{r} g_s((\bar{\bc}^{(l)})^T_s x) 
	- 
	\sum_{s=1}^r \sum_{k=1}^K a_k \cdot \sigma\left(\sum_{j=1}^d (b^{(l)})^{(0)}_{(s-1)\cdot K + k, j} \cdot x^{j} + 	(b^{(l)})^{(0)}_{(s-1)\cdot K + k, 0}) \right) - a_0|
	\\
	&&
	\leq
	\sum_{s=1}^{r}
	|g_s((\bar{\bc}^{(l)})^T_s x) 
	- 
	\sum_{k=1}^K a_k \cdot \sigma\left(\rho_n\cdot((\bar{\bc}^{(l)})^T_s x - (b^{(l)})_k\right) - a_0|
	\\
	&&
	\leq
	r \cdot 3 \cdot C \cdot \frac{(4 \cdot A \cdot \sqrt{d})^p}{(K-1)^p} + C \cdot (4 \cdot A \cdot \sqrt{d})^p \cdot (K-1)^{1-p} \cdot e^{- \frac{\rho_n \cdot A \cdot \sqrt{d}}{(n+1) \cdot (K-1)}}
	\\
	&&
	\leq
	const \cdot r \cdot C \cdot \frac{1}{K^p}
\end{eqnarray*}
Together this implies 
\begin{eqnarray*}
&&
\EXP \int |m_n(x)-m(x)|^2 \PROB_X (dx)
\\
&&
\leq
c_{11} \cdot \frac{(\log n)^3 \cdot K \cdot r}{n}
+ c_{12} \cdot r^2 \cdot C^2 \cdot \frac{1}{K^{2p}}
+
c_{13} \cdot
\EXP
\left\{
\min_{l=1,\dots,I_n}
\max_{s=1, \dots,r}
\|\bc_s - \bar{\bc}_s^{(l)} \|_{\infty}^{2p}
\right\}.
\end{eqnarray*}
The definition of $K$ implies
\[
c_{11} \cdot \frac{(\log n)^3 \cdot K \cdot r}{n}
+ c_{12} \cdot r^2 \cdot C^2 \cdot \frac{1}{K^{2p}}
\leq
c_{14}
\cdot
\left(
\frac{(\log n)^3}{n}
\right)^{\frac{2p}{2p+1}},
\]
hence it remains to show that we also have
\[
\EXP
\left\{
\min_{l=1,\dots,I_n}
\max_{s=1, \dots,r}
\|\bc_s - \bar{\bc}_s^{(l)}\|_{\infty}^{2p}
\right\}
\leq
c_{15}
\cdot
\left(
\frac{(\log n)^3}{n}
\right)^{\frac{2p}{2p+1}}
.
\]
By the random choice of the 
$\bar{\bc}_s^{(l)}$
we know for any $t \in (0,1]$
\begin{eqnarray*}
	\PROB\left\{
	\min_{l=1,\dots,I_n}
	\max_{s=1, \dots,r}
	\|\bc_s - \bar{\bc}_s^{(l)}\|_{\infty}
	>t
	\right\}
	&=&
	\prod_{i=1}^{I_n} \left(1-\PROB\left\{
	\max_{s=1, \dots,r}
	\|\bc_s - \bar{\bc}_s^{(i)}\|_{\infty}
	\leq t
	\right\} \right))
	\\
	&\leq&
	\left(
	1-t^{r \cdot d}
	\right)^{I_n}
\end{eqnarray*}
from which we conclude
\begin{eqnarray*}
	&&
  \EXP
  \left\{
  \min_{l=1,\dots,I_n}
  \max_{s=1, \dots,r}
  \|\bc_s - \bar{\bc}_s^{(l)}\|_{\infty}^{2p}
  \right\}
  \\
  &&
  =
  \int_{0}^{1} \PROB \left(
  \min_{l=1,\dots,I_n}
  \max_{s=1, \dots,r}
  \|\bc_s - \bar{\bc}_s^{(l)}\|_{\infty}^{2p}
  > t
   \right)
  dt
  \\
  &&
  =
  \int_{0}^{1} \PROB \left(
  \min_{l=1,\dots,I_n}
  \max_{s=1, \dots,r}
  \|\bc_s - \bar{\bc}_s^{(l)}\|_{\infty}
  > t^{\frac{1}{2p}}
  \right)
  dt
  \\
  &&
  \leq
   \int_{0}^{1} \exp \left( - I_n \cdot t^{\frac{r\cdot d}{2p}} \right) dt
  \\
  &&
  \leq
  \frac{2p}{r \cdot d} \cdot I_n^{-\frac{2p}{r \cdot d}} \cdot  \int_{0}^{\infty} e^{-s} \cdot s^{\frac{2p}{r \cdot d} -1} ds  
  \\
  &&
  =
   \frac{2p}{r \cdot d} \cdot I_n^{-\frac{2p}{r \cdot d}} \cdot \Gamma \left(\frac{2p}{r \cdot d}\right)
  \\
  &&
  \leq
  c_{15}
  \cdot
  \left(
  \frac{(\log n)^3}{n}
  \right)^{\frac{2p}{2p+1}},
\end{eqnarray*}
where the last inequatlity holds by assumption, since $p,r,d > 0$ are fixed.
Summarizing the above results we get the assertion.
\hfill $\Box$

\newpage

\begin{appendix}
  
\section{Supplementary material}

\subsection{Proof of Lemma \ref{le1}.}
           {\bf a)}
    For $s \in [0,1]$ set
    \[
H(s)=F(\ba^{(t)} + s \cdot (\ba^{(t+1)}-\ba^{(t)})).
    \]
    Then the fundamental theorem of calculus, the chain rule,
    the Cauchy-Schwarz inequality and assumption (\ref{le1eq2})
    imply
    \begin{eqnarray*}
      &&F(\ba^{(t+1)})-F(\ba^{(t)}) = H(1)-H(0) 
      = \int_0^1 H^\prime (s) \, ds \\
      &&= \int_0^1 (\nabla_\ba F)(\ba^{(t)} + s \cdot (\ba^{(t+1)}-\ba^{(t)}))
      \cdot (\ba^{(t+1)}-\ba^{(t)}) \, ds\\
      &&=\int_0^1 \left(
      (\nabla_\ba F)(\ba^{(t)} + s \cdot (\ba^{(t+1)}-\ba^{(t)}))
      -
      (\nabla_\ba F)(\ba^{(t)})
      \right)
      \cdot (\ba^{(t+1)}-\ba^{(t)}) \, ds\\
      && \quad
      +
     \int_0^1 
      (\nabla_\ba F)(\ba^{(t)})
      \cdot (\ba^{(t+1)}-\ba^{(t)}) \, ds
      \\
      &&
      \leq
      \int_0^1 L_n \cdot \|s \cdot (\ba^{(t+1)}-\ba^{(t)}) \| \cdot
      \|\ba^{(t+1)}-\ba^{(t)}\| \, ds\\
      && \quad +
      (\nabla_\ba F)(\ba^{(t)}) \cdot
      (\ba^{(t+1)}-\ba^{(t)})
      \\
      &&
      =
      \frac{L_n}{2} \cdot \|\ba^{(t+1)}-\ba^{(t)}\|^2
      +
      (\nabla_\ba F)(\ba^{(t)}) \cdot
      (\ba^{(t+1)}-\ba^{(t)}).
    \end{eqnarray*}
    Using (\ref{se5eq2}) and (\ref{le1eq1}) we get
    \begin{eqnarray*}
    F(\ba^{(t+1)})-F(\ba^{(t)})
    &\leq&
    \frac{L_n}{2} \cdot \lambda_n^2 \cdot \|  (\nabla_\ba F)(\ba^{(t)})\|^2
    - \lambda_n \|  (\nabla_\ba F)(\ba^{(t)})\|^2\\
    &
    =&
- \frac{1}{2 \cdot L_n}
      \cdot \|(\nabla_\ba F)(\ba^{(t)})\|^2.    
    \end{eqnarray*}

    \noindent
        {\bf b)}
        From {\bf a)} and (\ref{le1eq3}) we get
        \begin{eqnarray*}
          &&
                F(\ba^{(t+1)})-F(\ba_{opt})
                \\
                && \leq
                F(\ba^{(t)}) -F(\ba_{opt}) - \frac{1}{2 \cdot L_n}
                \cdot \|(\nabla_\ba F)(\ba^{(t)})\|^2
                \\
                && \leq
                                F(\ba^{(t)}) -F(\ba_{opt}) - \frac{1}{2 \cdot L_n}
                                \cdot
                                \rho_n \cdot (F(\ba^{(t)})-F(\ba_{opt}))\\
                                &&
                                =
                                      \left(1 - \frac{\rho_n}{2 \cdot L_n} \right)
      \cdot
(      F(\ba^{(t)})-F(\ba_{opt})
).      
        \end{eqnarray*}
        \quad
        \hfill $\Box$

\subsection{Proof of Lemma \ref{le9}}
In the proof we use the following error decomposition:
\beq
        \int |m_n(x) - m(x)|^2 \PROB_X (dx)&& \qquad \qquad\qquad \qquad\qquad \qquad\qquad \qquad\qquad \qquad
\eeq
\vspace{-8mm}
\beq
        &=&\Big[ \textbf{E}\Big\{|m_n(X) -Y|^2|\mathcal{D}_n\Big\} - \textbf{E} \Big\{|m(X)-Y|^2 \Big\} \\
  &&\qquad -
\Big(
  \textbf{E} \Big\{|m_n(X) - T_{\beta_n} Y|^2|\mathcal{D}_n \Big\} - \textbf{E} \Big\{|m_{\beta_n}(X)-
        T_{\beta_n} Y|^2 \Big\} \Big) \Big]\\
        &&+\Bigg[ \textbf{E} \Big\{|m_n(X) - T_{\beta_n} Y|^2|\mathcal{D}_n \Big\} - \textbf{E} \Big\{|m_{\beta_n}(X)
        - T_{\beta_n} Y|^2 \Big\}\\
        &&\qquad  - 2 \cdot \frac{1}{n} \sum_{i=1}^n \Big(|m_n(X_i) - T_{\beta_n} Y_i|^2 - |m_{\beta_n} (X_i)
  - T_{\beta_n} Y_i|^2 \Big) \Bigg]\\
        && + \left[ 2 \cdot \frac{1}{n} \sum_{i=1}^n |m_n(X_i) - T_{\beta_n} Y_i|^2
                                                - 2\cdot \frac{1}{n} \sum_{i=1}^n |m_{\beta_n} (X_i) - T_{\beta_n} Y_i|^2 \right. \\
        && \qquad \left. -\left( 2\cdot \frac{1}{n} \sum_{i=1}^n |m_n(X_i)-Y_i|^2 - 2\cdot \frac{1}{n} \sum_{i=1}^n|m(X_i)-Y_i|^2 \right)\right]\\
        && + \left[2\left( \frac{1}{n} \sum_{i=1}^n |m_n(X_i)-Y_i|^2 -\frac{1}{n} \sum_{i=1}^n|m(X_i)-Y_i|^2\right)\right]\\
        &=& \sum_{i=1}^4 T_{i,n},
\eeq
where $T_{\beta_n} Y$ is the truncated version of $Y$ and $m_{\beta_n}$ is the regression function of $T_{\beta_n} Y$, i.e.,
\[ m_{\beta_n}(x) = \textbf {E} \Big\{ T_{\beta_n} Y|X=x\Big\}.\]
We start with bounding $T_{1,n}$. By using $a^2-b^2 =(a-b)(a+b)$ we get
\beq
        T_{1,n}&=& \textbf{E} \Big\{|m_n(X) -Y|^2 - |m_n(X)- T_{\beta_n} Y|^2 \Big|\mathcal{D}_n \Big \}\\
        && -\textbf{E}\Big\{|m(X) - Y|^2 - |m_{\beta_n}(X)- T_{\beta_n} Y|^2 \Big\}\\
        &=& \textbf{E}\Big\{ (T_{\beta_n} Y -Y)  (2m_n(X)-Y-T_{\beta_n} Y) \Big|\mathcal D_n \Big \}\\
        && - \textbf{E}\Big\{ \Big( (m(X)-m_{\beta_n}(X)) + (T_{\beta_n} Y -Y) \Big) \Big( m(X) + m_{\beta_n}(X) -Y -T_{\beta_n} Y\Big) \Big\}\\
        &=& T_{5,n} + T_{6,n}.
\eeq
With the Cauchy-Schwarz inequality and
\beqm \label{expdurchexp}
I_{\{|Y|>\beta_n\}} \leq \frac{\exp (c_{2}/2 \cdot |Y|^2)}{\exp(c_{2}/2 \cdot \beta_n^2)}
\eeqm
 we conclude
\beq
|T_{5,n}| &\leq& \sqrt{\textbf{E}\big\{|T_{\beta_n} Y -Y|^2 \big\}} \cdot \sqrt{\textbf{E}\big\{|2m_n(X) -Y -T_{\beta_n} Y|^2\big|\mathcal D_n\big\}}\\
        &\leq & \sqrt{ \textbf{E}\big\{ |Y|^2 \cdot I_{\{|Y|>\beta_n\}} \big\}} \cdot \sqrt{\textbf{E}\big\{2\cdot |2m_n(X)
                -T_{\beta_n} Y|^2 + 2\cdot|Y|^2\big|\mathcal D_n\big\}}\\
        &\leq& \sqrt{ \textbf{E}\Bigg\{ |Y|^2 \cdot \frac{\exp(c_{2}/2\cdot|Y|^2)}{\exp(c_{2}/2\cdot \beta_n^2)} \Bigg\}} \\
        && \qquad \cdot
                 \sqrt{\textbf{E}\big\{2\cdot |2m_n(X) -T_{\beta_n} Y|^2 \big|\mathcal D_n \big\}+ 2\textbf{E}\big\{|Y|^2\big\}}\\
                 &\leq& \sqrt{\textbf{E}\Big\{ |Y|^2 \cdot \exp(c_{2}/2\cdot|Y|^2)\Big\}} \cdot
                \exp\left(-\frac{c_{2}\cdot \beta_n^2}{4}\right)
                  \cdot \sqrt{2 (3\beta_n)^2 + 2\textbf{E}\big\{|Y|^2\big\}}.
\eeq
With $ x \leq \exp(x)$ for $x \in \mathbb R$ we get
\[ |Y|^2 \leq \frac{2}{c_{2}} \cdot \exp \left( \frac{c_{2}}{2}\cdot |Y|^2\right) \]
and hence $\textbf{E}\Big\{ |Y|^2 \cdot \exp(c_{2}/2\cdot|Y|^2)\Big\}$ is bounded by
\[
\textbf{E} \left(  \frac{2}{c_{2}} \cdot \exp \left( c_{2}/2 \cdot |Y|^2\right)\cdot \exp(c_{2}/2 \cdot |Y|^2)\right) \leq
\textbf{E} \left(  \frac{2}{c_{2}} \cdot \exp \left( c_{2} \cdot |Y|^2 \right) \right) \leq c_{16}
\]
which is less than infinity by the assumptions of the lemma. Furthermore the third term is bounded by $\sqrt{18 \beta_n^2 + c_{17}}$  because
\beqm \label{Y^2}
\mathbf E(|Y|^2)\leq \mathbf E(1/c_{2} \cdot \exp( c_{2} \cdot |Y|^2) \leq c_{18} < \infty,
\eeqm
 which follows again as above. With the setting $\beta_n = c_{3} \cdot \log(n)  $ it follows for some constants $c_{19}, c_{20} >0$ that
\beq
        |T_{5,n}| &\leq& \sqrt{c_{16}} \cdot \exp \left( -c_{19} \cdot
          \log(n)^2\right) \cdot\sqrt{ (18\cdot c_{3}^2 \cdot (\log n)^2 +c_{17})}
	\leq c_{20} \cdot  \frac{\log(n)}{n}.
\eeq
From the Cauchy-Schwarz inequality we get
\beq
T_{6,n} &\leq& \sqrt{2\cdot \textbf{E}\Bigg\{ | (m(X)-m_{\beta_n}(X))|^2 \Big\}+ 2\cdot  \textbf{E} \Big\{ | (T_{\beta_n} Y -Y) |^2 \Bigg\}}\\
                                && \quad
                         \cdot \sqrt{\textbf {E} \Bigg\{ \Big| m(X) + m_{\beta_n}(X) -Y -T_{\beta_n} Y\Big|^2 \Bigg\}},
\eeq
where we can bound the second factor on the right-hand side in the above inequality in the same way we have bounded the second factor from $T_{5,n}$, because by assumption $||m||_\infty$ is bounded and furthermore $m_{\beta_n}$ is bounded by $\beta_n$. Thus we get for some constant $c_{21}>0$
 \[
 \sqrt{\textbf {E} \Bigg\{ \Big| m(X) + m_{\beta_n}(X) -Y -T_{\beta_n} Y\Big|^2 \Bigg\}} \leq c_{21} \cdot \log(n).
 \]
 Next we consider the first term. With Jensen's inequality it follows that
 \beq
 \textbf{E} \Big\{|m(X)-m_{\beta_n}(X)|^2 \Big\}
 &\leq& \textbf{E} \left\{ \textbf{E} \Big( |Y-T_{\beta_n} Y|^2 \Big| X \Big) \right\} = \textbf{E} \Big\{ |Y-T_{\beta_n} Y|^2 \Big\}.
 \eeq
 Hence we get
 \beq
 T_{6,n} &\leq& \sqrt{ 4 \cdot \textbf{E} \left\{ |Y-T_{\beta_n} Y|^2 \right\}} \cdot c_{21} \cdot \log(n)
 \eeq
 and therefore with the calculations from $T_{5,n}$ it follows that
 $T_{6,n}\leq c_{23} \cdot \log(n)/n$
 for some constant $c_{23}>0$. Altogether we get
 \[T_{1,n} \leq c_{24} \cdot \frac{\log(n)}{n}\]
 for some constant $c_{24}>0$.\\
Next we consider $T_{2,n}$ and conclude for $t>0$
\begin{align*}
\mathbf P\{T_{2,n}>t\}
        &\leq\mathbf P\left\{ \exists f\in T_{\beta_n} \mathcal{F}_n:
        \textbf{E} \left(\left|\frac{f(X)}{\beta_n}-\frac{T_{\beta_n} Y}{\beta_n}\right|^2\right)
        - \textbf{E} \left(\left|\frac{m_{\beta_n}(X)}{\beta_n}-\frac{T_{\beta_n} Y}{\beta_n} \right|^2\right) \right.\\
        & \qquad  -\frac{1}{n} \sum_{i=1}^n \left( \left| \frac{f(X_i)}{\beta_n} - \frac{T_{\beta_n} Y_i}{\beta_n}\right|^2
        -\left|\frac{m_{\beta_n}(X_i)}{\beta_n} - \frac{T_{\beta_n} Y_i}{\beta_n} \right|^2 \right)\\
        & \qquad > \frac{1}{2} \left( \frac{t}{\beta_n^2}
         +\textbf{E} \left(\left|\frac{f(X)}{\beta_n} -\frac{T_{\beta_n} Y}{\beta_n}\right|^2\right) \left.- \textbf{E} \left(\left| \frac{m_{\beta_n}(X)}{\beta_n}
                -\frac{ T_{\beta_n} Y}{\beta_n} \right|^2\right)\right)\right\},
\end{align*}
where $T_{\beta_n} \mathcal{F}_n$ is defined as $\left\{T_{\beta_n}f:f\in \F_n\right\}$.
Theorem 11.4 in Gy\"orfi et al. (2002) and the relation
\[
\mathcal N_1 \left( \delta , \left\{ \frac{1}{\beta_n} g :g \in \mathcal
    G \right\} , 
x_1^n \right)
\leq \mathcal N_1 \left( \delta \cdot \beta_n, \mathcal G, 
x_1^n \right)
\]
for an arbitrary function space $\G$ and $\delta >0$ lead to
\begin{align*}
\PROB \{T_{2,n}>t\}
        &\leq  14\cdot  \sup_{x_1^n}
\mathcal{N}_1 \left(\frac{t}{80\cdot\beta_n},
          \mathcal{F}_n, 
x_1^n
\right)
                \cdot \exp \left( - \frac{n}{5136 \cdot \beta_n^2} \cdot t \right).
\end{align*}
 Since the covering number is decreasing in $t$, we can conclude for $\varepsilon_n\geq \frac{80}{n}$
\begin{align*}
\EXP (T_{2,n}) &\leq \varepsilon_n + \int_{\varepsilon_n}^\infty \PROB \{T_{2,n}>t\} dt\\
&\leq \varepsilon_n + 14\cdot 
\sup_{x_1^n}
\mathcal{N}_1 \left(\frac{1}{n\cdot\beta_n}, \mathcal{F}_n,
  x_1^n
\right)
                \cdot \exp \left( - \frac{n}{5136 \cdot \beta_n^2} \cdot \varepsilon_n \right)\cdot \frac{5136 \cdot \beta_n^2}{n}.
\end{align*}
Choosing
\[
\varepsilon_n=\frac{5136\cdot \beta_n^2}{n}\cdot\log\left(
14
\cdot  
\sup_{x_1^n}
\mathcal{N}_1 \left(\frac{1}{n\cdot\beta_n}, \mathcal{F}_n,
  x_1^n
\right)\right)
\]
(which satisfies the necessary condition
$\varepsilon_n\geq\frac{80}{n}$ if the constant $c_{3}$ in the
definition of $\beta_n$ is not too small)
minimizes the right-hand side and implies
\begin{align*}
  \EXP (T_{2,n}) &\leq \frac{c_{25}\cdot \log(n)^2\cdot
    \log\left(
\sup_{x_1^n}
    \mathcal{N}_1 \left(\frac{1}{n\cdot\beta_n},
\mathcal{F}_n, x_1^n\right)
\right)
  }{n}.
\end{align*}
\noindent By bounding $T_{3,n}$ similarly to $T_{1,n}$ we get
\beq
\mathbf E (T_{3,n}) & \leq & c_{26} \cdot \frac{\log(n)}{n}
\eeq
for some large enough constant $c_{26}>0$ and hence we get in total
\beq
\mathbf E \left( \sum_{i=1}^3 T_{i,n} \right) &\leq&  \frac{c_{27} \cdot (\log n)^2
  \cdot
  \left(
\log\left(
\sup_{x_1^n} \mathcal{N}_1 \left(\frac{1}{n\cdot\beta_n},
\mathcal{F}_n, x_1^n \right)
\right)
+1
\right)
}{n}
\eeq
for some sufficient large constant $c_{27}>0$.

We finish the proof by bounding $T_{4,n}$.
Let $A_n$ be the event, that there exists $i \in \{1,...,n\}$ such that $|Y_i|>\beta_n$ and let $I_{A_n}$ be the indicator function of $A_n$. Then we get
\beq
        \textbf{E} (T_{4,n}) &\leq& 2 \cdot\mathbf E \left( \frac{1}{n} \sum_{i=1}^n |m_n(X_i)-Y_i|^2 \cdot I_{A_n}\right)\\
    && + 2 \cdot \mathbf E  \left( \frac{1}{n} \sum_{i=1}^n |m_n(X_i)-Y_i|^2 \cdot  I_{A_n^c}
        - \frac{1}{n} \sum_{i=1}^n |m(X_i)-Y_i|^2  \right)\\
        &=& 2 \cdot\mathbf E \left( |m_n(X_1)-Y_1|^2 \cdot I_{A_n}\right)\\
    && + 2 \cdot \mathbf E  \left( \frac{1}{n} \sum_{i=1}^n |m_n(X_i)-Y_i|^2 \cdot  I_{A_n^c}
        - \frac{1}{n} \sum_{i=1}^n |m(X_i)-Y_i|^2  \right)\\
       &=& T_{7,n} + T_{8,n}.
\eeq
With the Cauchy-Schwarz inequality we get for $T_{7,n}$
\beq
\frac{1}{2} \cdot T_{7,n}
&\leq& \sqrt{ \mathbf E \left( \left( |m_n(X_1) - Y_1|^2\right)^2 \right)} \cdot  \sqrt{ \mathbf P(A_n)}\\
&\leq& \sqrt{ \mathbf E \left( \left( 2 |m_n(X_1)|^2  + 2 |Y_1|^2\right)^2 \right) } \cdot \sqrt{ n \cdot \mathbf P \{|Y_1| >\beta_n \} }\\
&\leq & \sqrt{ \mathbf E \left(  8 |m_n(X_1)|^4  + 8 |Y_1|^4 \right) } \cdot \sqrt{ n \cdot
                    \frac{\mathbf E \left( \exp(c_{2}\cdot |Y_1|^2)\right)}{\exp( c_{2} \cdot \beta_n^2)} },
\eeq
where the last inequality follows as in the proof of
inequality (\ref{expdurchexp}).
With $x \leq \exp(x)$ for $x \in \mathbb R$ we get
\beq
\mathbf E \left( |Y|^4 \right) &=&\mathbf E \left( |Y|^2 \cdot |Y|^2\right)
                                                        \leq \mathbf E \left( \frac{2}{c_{2}} \cdot \exp\left( \frac{c_{2}}{2} \cdot |Y|^2\right)
                                                         \cdot \frac{2}{c_{2}} \cdot \exp\left( \frac{c_{2}}{2} \cdot |Y|^2\right) \right)\\
                                                         &=& \frac{4}{c_{2}^2} \cdot \mathbf E \left( \exp\left( c_{2} \cdot |Y|^2\right)\right),
\eeq
which is less than infinity by assumption (\ref{subgaus}) of the lemma.
Furthermore $||m_n||_\infty$ is bounded by $\beta_n$ and therefore the first factor is bounded by
\[
c_{28} \cdot \beta_n^2 = c_{29} \cdot (\log n)^2
\]
for some constant $c_{29}>0$. The second factor is bounded by $1/n$, because by the assumptions of the lemma $\mathbf E\left( \exp\left( c_{2}\cdot |Y_1|^2\right) \right)$ is bounded by some constant $c_{30} <\infty$ and hence
 we get
\beq
\sqrt{ n \cdot \frac{\mathbf E \left( \exp(c_{2}\cdot |Y_1|^2)\right)}{\exp( c_{2} \cdot \beta_n^2)} } &\leq&
        \sqrt{n} \cdot \frac{ \sqrt{c_{30}}}{\sqrt{\exp(c_{2} \cdot \beta_n^2)}}
        \leq \frac{\sqrt{n} \cdot
\sqrt{c_{30}}}{\exp((c_{2} \cdot c_{3}^2 \cdot (\log n)^2)/2)}.
\eeq
Since $\exp( - c \cdot \log(n)^2) = O(n^{-2})$ for any $c>0$, we get altogether
\beq
 T_{7,n}  &\leq& c_{31} \cdot \frac{(\log n)^2 \sqrt{n}}{n^2} \leq
 c_{32} \cdot 
\frac{(\log n)^2}{n}.
\eeq
With the definition of $A_n^c$ and $\tilde m_n$ defined as in the assumptions of this lemma we conclude
\beq
T_{8,n} &\leq& 2 \cdot \mathbf E  \left( \frac{1}{n} \sum_{i=1}^n |\tilde{m}_n(X_i)-Y_i|^2 \cdot  I_{A_n^c}
        - \frac{1}{n} \sum_{i=1}^n |m(X_i)-Y_i|^2  \right)\\
&\leq & 2 \cdot\mathbf E \left(
\min_{l \in \Theta_n}
\frac{1}{n} \sum_{i=1}^n | g_{n,l} (X_i)-Y_i|^2 + pen_n(g_{n,l})
+ {\epsilon}_{n,l}
        - \frac{1}{n} \sum_{i=1}^n |m(X_i)-Y_i|^2
\right)
\eeq
because $|T_\beta z - y| \leq |z-y|$ holds for $|y|\leq \beta$.
Since $pen_n(g_{n,l})$ and $\epsilon_{n,l}$ are deterministic terms
and
and since $g_{n,l}$ are independent
of $Y_1$, \dots, $Y_n$ given $\mathcal{D}_{n,r}$
we get that
\begin{eqnarray*}
	&&
	\mathbf E \left(
	\min_{l \in \Theta_n}
	\frac{1}{n} \sum_{i=1}^n | g_{n,l} (X_i)-Y_i|^2 + pen_n(g_{n,l})
	+ {\epsilon}_{n,l}
	- \frac{1}{n} \sum_{i=1}^n |m(X_i)-Y_i|^2
	\right)
	\\
	&&
	=
	\EXP \left(
	\mathbf E \left(
	\min_{l \in \Theta_n}
	\frac{1}{n} \sum_{i=1}^n | g_{n,l} (X_i)-Y_i|^2 
	- \frac{1}{n} \sum_{i=1}^n |m(X_i)-Y_i|^2
	+ pen_n(g_{n,l})
	+ {\epsilon}_{n,l}
	\mid
	\mathcal{D}_{n,r}
	\right)
	\right)
	\\
	&&
	\leq
	\EXP \left(
	\min_{l \in \Theta_n}
	\mathbf E \left(
	\frac{1}{n} \sum_{i=1}^n | g_{n,l} (X_i)-Y_i|^2 
	- \frac{1}{n} \sum_{i=1}^n |m(X_i)-Y_i|^2
	\mid
	\mathcal{D}_{n,r}
	\right)
	+ pen_n(g_{n,l})
	+ {\epsilon}_{n,l}
	\right)
	\\
	&&
	=
	\EXP \Bigg(
	\min_{l \in \Theta_n}
	\mathbf E \left(
	\frac{1}{n} \sum_{i=1}^n | g_{n,l} (X_i)-Y_i|^2 
	\mid
	\mathcal{D}_{n,r}
	\right)
	- 
	\EXP \left(
	\frac{1}{n} \sum_{i=1}^n |m(X_i)-Y_i|^2
	\mid
	\mathcal{D}_{n,r}
	\right)
	\\
	&&
	\qquad
	+ pen_n(g_{n,l})
	+ {\epsilon}_{n,l}
	\Bigg)
	\\
	&&
	=
	\EXP \Bigg(
	\min_{l \in \Theta_n}
	\mathbf E \left(
	\frac{1}{n} \sum_{i=1}^n | (g_{n,l} (X_i) - m(X_i)) + (m(X_i) -Y_i)|^2 
	\mid
	\mathcal{D}_{n,r}
	\right)
	\\
	&&
	\qquad
	-
	\EXP \left(
	\frac{1}{n} \sum_{i=1}^n |m(X_i)-Y_i|^2
	\mid
	\mathcal{D}_{n,r}
	\right)
	+ pen_n(g_{n,l})
	+ {\epsilon}_{n,l}
	\Bigg)
	\\
	&&
	=
	\EXP \Bigg(
	\min_{l \in \Theta_n}
	\mathbf E \left(
	\frac{1}{n} \sum_{i=1}^n |g_{n,l} (X_i) - m(X_i)|^2 
	\mid
	\mathcal{D}_{n,r}
	\right)
	+
	\mathbf E \left(
	\frac{1}{n} \sum_{i=1}^n  |m(X_i) -Y_i|^2 
	\mid
	\mathcal{D}_{n,r}
	\right)
	\\
	&&
	\qquad
	-
	\EXP \left(
	\frac{1}{n} \sum_{i=1}^n |m(X_i)-Y_i|^2
	\mid
	\mathcal{D}_{n,r}
	\right)
	+ pen_n(g_{n,l})
	+ {\epsilon}_{n,l}
	\Bigg)
	\\
	&&
	=
	\EXP \Bigg(
	\min_{l \in \Theta_n}
	\frac{1}{n} \sum_{i=1}^n |g_{n,l} (X_i) - m(X_i)|^2 
	+ pen_n(g_{n,l})
	+ {\epsilon}_{n,l}
	\Bigg)
\end{eqnarray*}
where the fourth equality holds
since the mixed term is
\begin{eqnarray*}
	&&
	\mathbf E \left(
	\frac{1}{n} \sum_{i=1}^n (g_{n,l} (X_i)-m(X_i)) \cdot (m(X_i)-Y_i) 
	\mid
	\mathcal{D}_{n,r}
	\right)
	\\
	&&
	=
	\mathbf E \left(
	\frac{1}{n} \sum_{i=1}^n (g_{n,l} (X_i)-m(X_i)) \cdot
	\mathbf E \left(
	(m(X_i)-Y_i) 
	\mid
	\mathcal{D}_{n,r}
	\right)
	\right)
	\\
	&&
	=
	\mathbf E \left(
	\frac{1}{n} \sum_{i=1}^n (g_{n,l} (X_i)-m(X_i)) \cdot
	\mathbf E \left(
	(m(X_i)-Y_i) 
	\mid
	X_i
	\right)
	\right)
	\\
	&&
	=
	0
\end{eqnarray*}
Hence,
\beq
&&
\mathbf E (T_{4,n}) \\
&&\leq c_{32} \cdot \frac{(\log n)^2}{n}
\\
&&
\quad
+
                        2 \cdot \mathbf E \left(  
\min_{l \in \Theta_n}
\frac{1}{n} \sum_{i=1}^n | g_{n,l}(X_i)-Y_i|^2 + pen_n(g_{n,l})
+ {\epsilon}_{n,l}
                        - \frac{1}{n} \sum_{i=1}^n |m(X_i)-Y_i|^2  \right)
                        \\
&&\leq 
c_{32} \cdot \frac{(\log n)^2}{n}
+
2 \cdot
	\EXP \left(
\min_{l \in \Theta_n}
\frac{1}{n} \sum_{i=1}^n |g_{n,l} (X_i) - m(X_i)|^2 
+ pen_n(g_{n,l})
+ {\epsilon}_{n,l}
\right)
\eeq
holds.
In combination with the other considerations in the proof this implies the assertion of Lemma \ref{le9}.\hfill $\Box$\\

\end{appendix}

\end{document}